\documentclass{article}

\usepackage{amsfonts}
\usepackage{amsmath}
\usepackage{amssymb}
\usepackage{amsthm}
\usepackage{enumerate}
\usepackage{psfrag}
\usepackage{txfonts}

\newtheorem{definition}{Definition}[section]
\newtheorem{theorem}{Theorem}[section]
\newtheorem{proposition}{Proposition}[section]

\newtheorem{corollary}{Corollary}[section]

\title{Minimum Number of Colors:\\
 the Turk's Head Knots Case Study}
\author{Pedro Lopes$\sp{1, 2}$ and Jo\~ao Matias$\sp{2}$\\
$\sp{1}$Center for Mathematical Analysis, Geometry and Dynamical Systems\\
$\sp{2}$Department of Mathematics\\
Instituto Superior T\'ecnico\\
Technical University of Lisbon\\
Av. Rovisco Pais\\
1049-001 Lisbon\\
Portugal\\
        \texttt{pelopes@math.ist.utl.pt  joao.matias@ist.utl.pt}\\}

\date{February 25, 2010}
\begin{document}

\maketitle

\begin{abstract}
The minimum number of colors is a challenging knot invariant since, by definition, its calculation requires taking the minimum
over infinitely many minima. In this article we estimate and in some cases calculate the minimum number of colors for the Turk's
head knots on three strands.
\end{abstract}

\bigbreak


\bigbreak


\bigbreak

\section{Introduction}

\noindent

\begin{figure}[!ht]
    \psfrag{o}{\huge $\text{over-arc}$}
    \psfrag{u}{\huge $\text{under-arcs}$}
    \psfrag{2}{\huge $2$}
    \psfrag{3}{\huge $3$}
    \psfrag{4}{\huge $4$}
    \psfrag{7}{\huge $7$}
    \psfrag{2-}{\huge A $2$-coloring of $THK(3, 3)$ with $2$ colors}
    \psfrag{5-}{\huge A $5$-coloring of $THK(2, 3)$ with $4$ colors}
    \centerline{\scalebox{.50}{\includegraphics{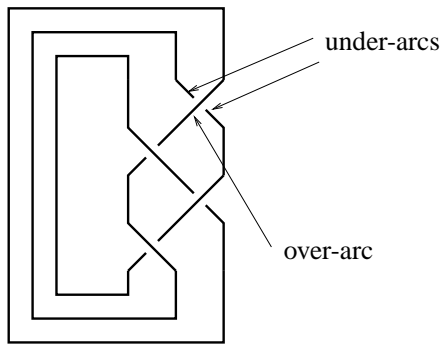}}}
    \caption{A knot diagram specifying over-arc and under-arcs at the indicated crossing}\label{fig:intro0}
\end{figure}

Knots (see \cite{lhKauffman} for further information) are embeddings of the circle into $3$-space modulo deformation  of the embedding. The classification of knots modulo deformations is still an open problem. Knots are usually represented via their projections onto the plane plus extra information: at the crossings in the projection, the line that goes under is broken to denote farther distance from the observer. The network of arcs so obtained is called {\it knot diagram} (see Figure \ref{fig:intro0}). The Reidemeister theorem establishes that if two knot diagrams are related by a finite number of the so-called Reidemeister moves, then the knots they stem from are deformable into each other (and conversely). Suppose a mathematical object is assigned to each knot diagram in such a way that it is (essentially) invariant under the Reidemeister moves. This mathematical object constitutes a {\it knot invariant} in the sense that if different objects are assigned to diagrams $D$ and $D'$, then the knots they stem from are not deformable into each other.

We now elaborate on Fox colorings since they give rise to the  the knot invariant we are concerned with in this article, the minimum number of colors.

\bigbreak

\begin{figure}[!ht]
    \psfrag{x}{\huge $x$}
    \psfrag{y}{\huge $y$}
    \psfrag{z}{\huge $z$}
    \psfrag{3}{\huge $3$}
    \psfrag{4}{\huge $4$}
    \psfrag{7}{\huge $7$}
    \psfrag{2-}{\huge A $2$-coloring of $THK(3, 3)$ with $2$ colors}
    \psfrag{5-}{\huge A $5$-coloring of $THK(2, 3)$ with $4$ colors}
    \centerline{\scalebox{.50}{\includegraphics{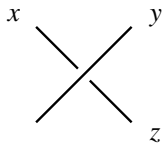}}}
    \caption{Arcs as variables $x, y$ and $z$; the crossing where they meet dictates the equation $2y-x-z=0$}\label{fig:cross}
\end{figure}

Given an integer $r>1$ and a knot diagram, $D\sb{K}$, of a given knot $K$, a (Fox) $r$-coloring of $D\sb{K}$ (\cite{Fox}) is an assignment of integers mod $r$ (i.e., numbers from $\mathbb{Z}\sb{r} = \{ 0, 1, 2, \dots , r-1  \}$ mod $r$) such that, at each crossing of $D\sb{K}$, the equality ``twice the color on the over-arc equals the sum of the colors on the under-arcs'' holds (mod $r$). The Fox $r$-colorings can be alternatively envisaged as follows. Assign a variable to each arc of $D\sb{K}$ and at each crossing read off the equation $2y-x-z=0$ where $y$ is the variable assigned to the over-arc and $x$ and $z$ are the variables assigned to the under-arcs (Figure \ref{fig:cross}).
A system of linear homogeneous equations is thus associated to each knot diagram. The solutions of this system of equations mod $r$ constitute the $r$-colorings of $D\sb{K}$. There are always the trivial colorings i.e., the colorings where each arc is assigned the same color. Upon performing a Reidemeister move on $D\sb{K}$ endowed with an $r$-coloring, we can consistently assign colors to the arcs on the transformed portion of diagram in a unique way so that we obtain an $r$-coloring on the new diagram. Furthermore, if we undo this Reidemeister move, we can reassign colors so that we obtain the original $r$-coloring on $D\sb{K}$. There is thus a bijection between the $r$-colorings of two knot diagrams related by a finite number of Reidemeister moves (\cite{pLopes}). Hence, the number of $r$-colorings is a knot invariant. Furthermore, this consistent assignment or reassignment of colors upon the performance of Reidemeister moves on colored diagrams takes trivial $r$-colorings to trivial $r$-colorings and non-trivial $r$-colorings to non-trivial $r$-colorings (non-trivial $r$-colorings being those that use at least two colors). In this way, the existence or not of non-trivial $r$-colorings is also a knot invariant.

\bigbreak

Assume, then, knot $K$ admits non-trivial $r$-colorings. What is the minimum number of colors needed to set up a non-trivial coloring over all diagrams of $K$?

\begin{definition} Let $r$ be an integer greater than $1$. Let $K$ be a knot, $D\sb{K}$ one of its diagrams. Assume $K$ admits
non-trivial $r$-colorings and let $n\sb{D\sb{K}, r}$ stand for the least number of colors it takes to set up a non-trivial
$r$-coloring of $D\sb{K}$. We call
\[
{\rm mincol}\sb{r}(K):=\min \, \{ \, n\sb{D\sb{K}, r} \, \, \big| \, \, D\sb{K} \text{ is a diagram of } K  \}
\]
the \emph{minimum number of colors of $K$ mod $r$}
\end{definition}

This is tautologically a knot invariant. It is striking that, by definition, one would have to consider each and every one of the
infinitely many diagrams of the knot under study in order to calculate this invariant, apparently making this task impossible to execute.

\bigbreak

We remarked before that  $r$-colorings of a knot $K$ can be regarded as the solutions of a system of linear homogeneous equations over $\mathbb{Z}\sb{r}$ read off from the diagram of $K$ under study. Upon performance of Reidemeister moves on this diagram, the matrix of the coefficients of the linear homogeneous system of equations undergoes {\it elementary transformations} as described in \cite{Lickorish} on page $50$. In this way, the equivalence class of the matrix of the coefficients modulo elementary transformations is a knot invariant and so is any minor of this matrix. The first minor of this matrix is called the determinant of the knot $K$, $\det K$.

\bigbreak

The topic of minimum number of colors was set forth in \cite{Frank} where the Kauffman-Harary conjecture was presented. Given a prime $p$, this conjecture states that a non-trivial $p$-coloring on a reduced diagram of an alternating knot of prime determinant $p$, assigns different colors to different arcs. At the time of the writing of this article, there is an alleged proof of this conjecture available (\cite{msolis}).

There are other articles on minimum number of colors addressing the actual minimum for a given $r$-coloring. Satoh (\cite{satoh}) shows that any non-trivial $5$-coloring can be realized with $4$ colors; Oshiro (\cite{Oshiro}) shows that any $7$-coloring can be realized with $4$ colors. On the other hand, Saito (\cite{Saito}) gives a condition for the minimum number of colors in a non-trivial $p$-coloring to be greater than $4$, for prime $p>7$.

\bigbreak

In \cite{kl}, the torus knots of type $(2, n)$, the $T(2, n)$'s, were investigated. A formula for the number of $r$-colorings was established for each $n$, thereby allowing one to realize for which pairs $(r, n)$ there are non-trivial $r$-colorings of $T(2, n)$. For these cases, and relying on the features of modular arithmetic, estimates and sometimes actual minima were presented for the minimum number of colors. Furthermore, a sequence of transformations on the standard diagrams of the $T(2, n)$'s were defined that further decreased the number of colors in infinitely many cases.

The present article is a sequel to \cite{kl}. As a matter of fact our original intention was to mimic the techniques of \cite{kl}, applying them to the class of the Turk's head knots on $3$ strands, $\{ THK(3, n)\}\sb{n\in \mathbb{Z}\sp{+}}$. This is the collection of knots given by the braid closure of $\big(\sigma\sb{2}\sigma\sb{1}\sp{-1}\big) \sp{n}\in B\sb{3}$:
\[
THK(3, n) = \text{ Braid Closure } \bigg[ \big(\sigma\sb{2}\sigma\sb{1}\sp{-1}\big) \sp{n} \bigg]
\]
\begin{figure}[!ht]
    \centerline{\scalebox{.50}{\includegraphics{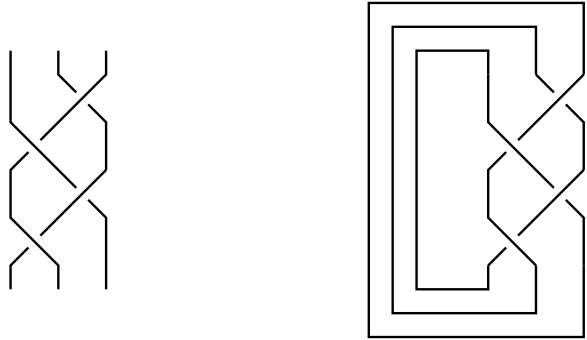}}}
    \caption{Illustrating the difference between braid (left-hand side) and the closure of the braid (right-hand side). The braid at issue is $\big(\sigma\sb{2}\sigma\sb{1}\sp{-1}\big) \sp{2}\in B\sb{3}$.}\label{fig:braidcl}
\end{figure}
See Figure \ref{fig:braidcl} for an illustrative example of the difference between ``braid'' and ``braid closure''.  As it happens, the increasing complexity of the Turk's head knots with respect to the torus knots of type $(2, n)$ turned the present work into anything but a mere application of the techniques of \cite{kl}.

See \cite{Birman} for further information on braids. We call the knot diagrams of $THK(3, n)$ given by the braid closure of $\big(\sigma\sb{2}\sigma\sb{1}\sp{-1}\big) \sp{n}\in B\sb{3}$, {\it standard diagrams} of $THK(3, n)$.

\bigbreak

There is another article which addresses colorings on Turk's head knots (\cite{THK}), but its goal is to show that families of  these knots satisfy the Kauffman-Harary's conjecture by showing their determinants are composite numbers. In this way, the current work addresses a different matter which is the calculation of minimum number of colors for different moduli.

\bigbreak

Here are the main results in this article. We remark that given two positive integers $a, b$, we let $(a, b)$ stand for their greatest common divisor.

\bigbreak

\begin{theorem}\label{thm:numbercols} Given positive integers $n$ and $r>1$, the number of $r$-colorings of $THK(3, n)$, denoted $\# col\sb{r} THK(3, n)$, is
\[
\# col\sb{r} THK(3, n) =
\begin{cases}
(u\sb{n-1}, r)\sp{2}r & \text{ if $n$ is odd}\\
(5u\sb{n-1}, r)(u\sb{n-1}, r)r & \text{ if $n$ is even}
\end{cases}
\]
where
\[
u\sb{n}=\frac{1}{\sqrt{5}}\Bigg[ \Bigg( \frac{1+\sqrt{5}}{2} \Bigg)\sp{n+2} - \Bigg( \frac{-1 + \sqrt{5}}{2} \Bigg)\sp{n}
 - \Bigg( \frac{1-\sqrt{5}}{2} \Bigg)\sp{n+2} + \Bigg( \frac{-1-\sqrt{5}}{2} \Bigg)\sp{n}\Bigg]
\]
\end{theorem}

\bigbreak

\begin{corollary}\label{cor:det}
\begin{equation*}
\det THK(3, n) \quad = \quad
\begin{cases}
(u\sb{n-1})\sp{2}& \quad \text{ if } n \text{ is odd}\\
5(u\sb{n-1})\sp{2}& \quad \text{ if } n \text{ is even}
\end{cases} \quad
\neq 0
\end{equation*}
\end{corollary}

\bigbreak

\begin{corollary}\label{cor:nontriv} There are non-trivial $r$-colorings of $THK(3, n)$ if, and only if,
\begin{itemize}
\item $(u\sb{n-1}, r)>1$ \qquad \qquad or

\item $n$ is even and $5\mid r$
\end{itemize}
\end{corollary}

\bigbreak

\begin{theorem}\label{thm:minexact}Let $n$ and $r$ be positive integers.
\begin{enumerate}
\item
\[
\qquad \qquad  2\mid r \text{  and } 3\mid n  \qquad \text{ if, and only if, }\qquad   mincol\sb{r}THK(3, n) = 2
\]
\item
 \[
 \Big( 2\nmid r  \text{  or }   3\nmid n\Big)\,   \text{  and }  \, 3\mid r \, \text{  and } \, 4\mid n \quad \text{ if, and only if, } \quad mincol\sb{r}THK(3, n) = 3
 \]
\item
\[
\Big( 2\nmid r  \text{ or }  3\nmid n \Big) \, \text{ and } \, \Big( 3\nmid r  \text{ or }  4\nmid n \Big)  \, \text{ and }  \,  \Big[\big( \, 5\mid r \text{ and } 2\mid n\, \big) \text{ or } \big(\,  7\mid r \, \text{ and }\,  8\mid n\, \big) \Big] \qquad
\]
\[
\text{ if, and only if,  }\qquad mincol\sb{r}THK(3, n) = 4
\]
\item If
\[
\Big( 2\nmid r \, \text{ or }  \, 3\nmid n \Big)  \, \text{ and }  \, \Big( 3\nmid r   \, \text{ or } \, 4\nmid n \Big)  \, \text{ and }  \, \Big( 5\nmid r  \, \text{  or }  \, 2\nmid n\Big)  \, \text {and }  \, \Big( 7\nmid r  \, \text{ or }  \, 8\nmid n\Big)
\]
\[
\text{and  } \qquad \Big( 11\mid r   \text{ and } 5\mid n \Big)  \qquad \text{ then } \qquad mincol\sb{r}THK(3, n) = 5
\]
\end{enumerate}
\end{theorem}

\bigbreak

It is easy to see that each statement in Theorem \ref{thm:minexact} gives rise to infinitely many knots for which the minimum number of colors modulo infinitely many $r$'s is exactly determined. For example, according to statement {\it 4.}, for any positive integers $m$ and $l$, we have:
\[
mincol\sb{\displaystyle{ 11\cdot 13\sp{l}}}THK(3, 5\cdot 17\sp{m}) = 5
\]

\bigbreak

\bigbreak

We introduce the mapping $\psi$ which associates to each color $r$, the least positive integer $n$ such that $r\mid u\sb{n-1}$.

\bigbreak

\begin{definition}\label{def:psi} For any integer $r>1$ set
\[
\psi (r) := \min \{ \, q\in \mathbb{Z}\sp{+}\, \big| \, r\mid u\sb{q-1}  \,   \}
\]
\end{definition}

\bigbreak

When $r$ is a prime other than $5$, $\psi (r)$ is the least number of $\sigma\sb{2}\sigma\sb{1}\sp{-1}$ that should be juxtaposed in order to obtain a non trivial $r$-coloring, namely in $THK(3, \psi (r))$.

\bigbreak

\begin{theorem}\label{thm:minestimates} Let $p$ be a prime greater than $11$.
\begin{enumerate}
\item Assume $\psi (p)$ is odd.
\begin{enumerate}
\item If $5\sp{\frac{p-1}{2}}=\sb{p}-1$, then $mincol\sb{p}THK(3, \psi (p)) \leq \frac{p+1}{2}$
\item If $5\sp{\frac{p-1}{2}}=\sb{p}1$, then $mincol\sb{p}THK(3, \psi (p)) \leq \frac{p-1}{2}$
\end{enumerate}
\item Assume $\psi (p)$ is even.
\begin{enumerate}
\item If $4\mid \psi (p)$, then $mincol\sb{p}THK(3, \psi (p)) \leq \psi ( p ) - 1$
\item If $4\nmid \psi (p)$, then $mincol\sb{p}THK(3, \psi (p)) \leq \psi ( p ) - 5$
\end{enumerate}
\end{enumerate}
\end{theorem}

\bigbreak

\begin{definition}\label{def:minpsi} For any positive integers $n$ and $r$, such that $(u\sb{n-1}, r)>1$, set
\[
\langle u\sb{n-1}, r \rangle\sb{\psi}
\]
to be the least common prime factor (greater than $5$) of $r$ and $u\sb{n-1}$, which minimizes $\psi$.
\end{definition}

\bigbreak

\begin{corollary}\label{cor:minestimatesgen} Let $r$ and $n$ be positive integers such that $(u\sb{n-1}, r)>1$. Set
\[
p:=\langle u\sb{n-1}, r \rangle\sb{\psi}
\]
Then,
\begin{enumerate}
\item Assume $\psi (p)$ is odd.
\begin{enumerate}
\item If $5\sp{\frac{p-1}{2}}=\sb{p}-1$, then $mincol\sb{r}THK(3, n) \leq \frac{p+1}{2}$
\item If $5\sp{\frac{p-1}{2}}=\sb{p}1$, then $mincol\sb{r}THK(3, n) \leq \frac{p-1}{2}$
\end{enumerate}
\item Assume $\psi (p)$ is even.
\begin{enumerate}
\item If $4\mid \psi (p)$, then $mincol\sb{r}THK(3, n) \leq \psi ( p ) - 1$
\item If $4\nmid \psi (p)$, then $mincol\sb{r}THK(3, n) \leq \psi ( p ) - 5$
\end{enumerate}
\end{enumerate}
\end{corollary}

\bigbreak

We establish below (Corollary \ref{cor:pmidup}) that, for prime $p\neq 5$, if $\psi (p)$ is odd then it is bounded above by $(p+1)/2$, whereas if $\psi (p)$ is even, it is bounded above by $p+1$. Then, for the odd $p$ case, Theorem \ref{thm:minestimates}
provides good estimates, roughly half the number of colors available ($p$). On the other hand, for the even $p$ case, the estimates are coarser. Above all for the $\psi (p) = p+1$ case where the results provide an estimate which equals the number of colors available ($p$).

We then ran a program in Mathematica to have  an idea of how many times this $\psi (p) = p+1$ situation occurs over all primes $p$. Our program did this for the first $100,000$ primes in steps of $10,000$. The results are displayed in Table \ref{Ta:ppsip}.
\begin{table}[h!!]
\begin{center}
\scalebox{.690}{\begin{tabular}{| c || c | c |  c | c | c | c |  c | c |  c | c | }
\hline
$N\sp{p}$ & $10,000$ &   $20,000$  & $30,000$ &   $40,000$ & $50,000$ &   $60,000$  & $70,000$ &   $80,000$  & $90,000$ &   $100,000$\\ \hline\hline
$N\sb{\psi (p)}$ & $3,969$ &   $7910$ &  $11,853$   &   $15,760$   &  $19,738$   &   $23,661$   &  $27,589$   &   $31,499$   &  $35,404$   &   $39,343$\\ \hline\hline
$N\sb{\psi (p)}/N\sp{p}$  & $0.3969$ &   $0.3955$ &  $0.3951$   &   $0.394$   &  $0.39476$   &   $0.39435$   &  $0.394129$   &   $0.393738$   &  $0.393378$   &   $0.39343$\\ \hline
\end{tabular}}
\caption{$N\sp{p}$ is the number of consecutive primes; $N\sb{\psi (p)}$ is the number of primes $p$ for which $\psi (p) = p+1$.}\label{Ta:ppsip}
\end{center}
\end{table}
This seems to indicate that around $40\%$ of the primes, $p$, lead to $\psi (p)=p+1$.

We  ran another program to realize what is the percentage of distinct colors used over the total number of colors for the colorings we used in this situation. These colorings were induced by introducing either colors $0, 1, 0$ or colors $1, 2, 0$ on the top of the standard diagram of the $THK(3, \psi (p))$ for each of these $p$'s such that $\psi (p) = p+1$ and $p>7$. For the first $100,000$ such cases we had a minimum of $69.2308 \%$ and a maximum of $75,0004 \%$ of colors used.

\bigbreak

In Section \ref{sect:proofs} we prove Theorems \ref{thm:numbercols}, \ref{thm:minexact}, and \ref{thm:minestimates}, along with their corollaries.

\subsection{Acknowledgements}\label{subsect:ackn}

\noindent

P.L. acknowledges support by the Funda\c{c}\~{a}o para a Ci\^{e}ncia e a Tecnologia (FCT /
Portugal). P.L. and J.M. acknowledge support by the Gulbenkian Foundation (Portugal) in connection with the 2008/2009 edition of the programme ``Novos Talentos em Matem\'atica''.

\bigbreak

\section{Proofs}\label{sect:proofs}

\noindent

In this Section we provide the proofs of the results stated in the Introduction.

\subsection{Proof of Theorem \ref{thm:numbercols}}\label{subsect:thmnumbercols}

\noindent

\begin{figure}[!ht]
    \psfrag{a}{\huge $a$}
    \psfrag{b}{\huge $b$}
    \psfrag{c}{\huge $c$}
    \psfrag{x}{\huge $x_1 = 2a-c$}
    \psfrag{y}{\huge $y_1=a$}
    \psfrag{z}{\huge $z_1 = 2c-b$}
    \centerline{\scalebox{.50}{\includegraphics{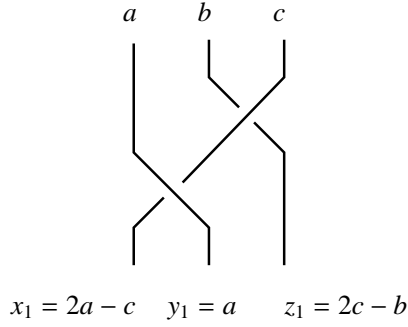}}}
    \caption{Propagation of colors down $\sigma\sb{2}\sigma\sb{1}\sp{-1}$}\label{fig:slq}
\end{figure}

In this Subsection we prove Theorem \ref{thm:numbercols}, which yields a formula for the number of $r$-colorings of $THK(3, n)$, along with its corollaries. For that, we start by studying the propagation of colors $a, b, c$ of an $r$-coloring, down $\big( \sigma\sb{2}\sigma\sb{1}\sp{-1}\big)\sp{n}$ (in Figure \ref{fig:slq} the case $n=1$ is displayed). We recall that this means that at each crossing, the equation ``twice the color on the over-arc equals the sum of colors on the under-arcs'' is satisfied modulo $r$. At each crossing, we use this rule to write the color of the lower under-arc in terms of the color of over-arc and of the color of the upper under-arc.

In Figure \ref{fig:slq} we illustrate the propagation of colors $a, b, c$ down $(\sigma\sb{2}\sigma\sb{1}\sp{-1})\sp{1}$, which, algebraically, translates into the system of equations \ref{eqn:algtrans}.

\begin{equation}\label{eqn:algtrans}
\begin{cases}
x_1= &2a-c\\
y_1= &a\\
z_1= &2c-b
\end{cases}\qquad\Longleftrightarrow\qquad
\left[
\begin{matrix}
x_1\\
y_1\\
z_1
\end{matrix}
\right]=
\left[
\begin{matrix}
2 &0 &-1\\
1 &0 &0\\
0 &-1 &2
\end{matrix}
\right]
\left[
\begin{matrix}
a\\
b\\
c
\end{matrix}
\right]
\end{equation}

In the sequel, we will use the following notation.

\begin{equation*}
C=
\left[
\begin{matrix}
2 &0 &-1\\
1 &0 &0\\
0 &-1 &2
\end{matrix}
\right]
\end{equation*}

We thus start by setting $x_{0}=a, y_{0}=b, z_{0}=c$ at the top of the braid, from left to right; we call this the {\it color input}. The colors $x_{n}, y_{n}, z_{n}$ (from left to right) after $\big( \sigma\sb{2}\sigma\sb{1}\sp{-1}\big)\sp{n}$ are:
\begin{equation*}
\left[
\begin{matrix}
x_{n}\\
y_{n}\\
z_{n}
\end{matrix}
\right]=C^{n}
\left[
\begin{matrix}
a\\
b\\
c
\end{matrix}
\right]
\end{equation*}

\bigbreak

\begin{proposition}\label{proposition:const}
Consider $THK(3, n)$ given by the braid closure of $(\sigma_{2}\sigma_{1}^{-1})^{n}$. Let $n$ be a non-negative integer. If the top strands of $(\sigma_{2}\sigma_{1}^{-1})^{n}$ are endowed with colors $a,b$ and $c$ (from left to right) then the following hold:

\begin{enumerate}[(i)]
 \item $x_{n}-y_{n}+z_{n}=a-b+c$.

\item $x_{n}=4x_{n-1}-4x_{n-2}+x_{n-3}$;\\
$y_{n}=4y_{n-1}-4y_{n-2}+y_{n-3}$;\\
$z_{n}=4z_{n-1}-4z_{n-2}+z_{n-3}$.

\item $x_{n}=3x_{n-1}-x_{n-2}-x_{0}+y_{0}-z_{0}$;\\
$y_{n}=3y_{n-1}-y_{n-2}-x_{0}+y_{0}-z_{0}$;\\
$z_{n}=3z_{n-1}-z_{n-2}-x_{0}+y_{0}-z_{0}$.
\end{enumerate}

\end{proposition}

\begin{proof}
To prove $(i)$ we note that:

\begin{equation*}
\left[
\begin{matrix}
1 &-1 &1
\end{matrix}
\right]
C=\left[
\begin{matrix}
1 &-1 &1
\end{matrix}
\right]
\left[
\begin{matrix}
2 &0 &-1\\
1 &0 &0\\
0 &-1 &2
\end{matrix}
\right]=
\left[
\begin{matrix}
1 &-1 &1
\end{matrix}
\right]
\end{equation*}

In this way, we obtain:

\begin{equation*}
\left[
\begin{matrix}
x_{n}-y_{n}+z_{n}
\end{matrix}
\right]
=
\left[
\begin{matrix}
1 &-1 &1
\end{matrix}
\right]
\left[
\begin{matrix}
x_{n}\\
y_{n}\\
z_{n}
\end{matrix}
\right]=
\left[
\begin{matrix}
1 &-1 &1
\end{matrix}
\right]
C^{m}
\left[
\begin{matrix}
a\\
b\\
c
\end{matrix}
\right]=
\left[
\begin{matrix}a-b+c
\end{matrix}
\right]
\end{equation*}

In order to prove $(ii)$ we note that, the characteristic polynomial of $C$ is
\[
det(C-I_{3}x)=-x^{3}+4x^{2}-4x+1
\]
Via the Cayley-Hamilton Theorem,

\begin{equation*}
(-C^{3}+4C^{2}-4C+I_{3})
\left[
\begin{matrix}
x_{n}\\
y_{n}\\
z_{n}
\end{matrix}
\right]=0
\end{equation*}
which amounts to
\begin{equation*}
-\left[
\begin{matrix}
x_{n+3}\\
y_{n+3}\\
z_{n+3}
\end{matrix}
\right]+4
\left[
\begin{matrix}
x_{n+2}\\
y_{n+2}\\
z_{n+2}
\end{matrix}
\right]-4
\left[
\begin{matrix}
x_{n+1}\\
y_{n+1}\\
z_{n+1}
\end
{matrix}
\right]+
\left[
\begin{matrix}
x_{n}\\
y_{n}\\
z_{n}
\end{matrix}
\right]=0
\end{equation*}

Thus yielding the indicated recurrence relations for $x\sb{n}, y\sb{n}, z\sb{n}$.

To prove $(iii)$ we start by using:
\begin{align*}
\left[
\begin{matrix}
x_{n}\\
y_{n}\\
z_{n}
\end{matrix}
\right]&=4
\left[
\begin{matrix}
x_{n-1}\\
y_{n-1}\\
z_{n-1}
\end{matrix}
\right]-4
\left[
\begin{matrix}
x_{n-2}\\
y_{n-2}\\
z_{n-2}
\end{matrix}
\right]+
\left[
\begin{matrix}
x_{n-3}\\
y_{n-3}\\
z_{n-3}
\end{matrix}
\right]=\\
&=3\left[
\begin{matrix}
x_{n-1}\\
y_{n-1}\\
z_{n-1}
\end{matrix}
\right]-
\left[
\begin{matrix}
x_{n-2}\\
y_{n-2}\\
z_{n-2}
\end{matrix}
\right]+
\left(
\left[
\begin{matrix}
2 &0 &-1\\
1 &0 &0\\
0 &-1 &2
\end{matrix}
\right]^{2}-3
\left[
\begin{matrix}
2 &0 &-1\\
1 &0 &0\\
0 &-1 &2
\end{matrix}
\right]+I_{3} \right)
\left[
\begin{matrix}
x_{n-3}\\
y_{n-3}\\
z_{n-3}
\end{matrix}
\right]=\\
&=3\left[
\begin{matrix}
x_{n-1}\\
y_{n-1}\\
z_{n-1}
\end{matrix}
\right]-
\left[
\begin{matrix}
x_{n-2}\\
y_{n-2}\\
z_{n-2}
\end{matrix}
\right]+
\left[
\begin{matrix}
 -1 &1 &-1\\
-1 &1 &-1\\
-1 &1 &-1
\end{matrix}
\right]
\left[
\begin{matrix}
 x_{n-3}\\
y_{n-3}\\
z_{n-3}
\end{matrix}
\right]= \\
&=3\left[
\begin{matrix}
x_{n-1}\\
y_{n-1}\\
z_{n-1}
\end{matrix}
\right]-
\left[
\begin{matrix}
x_{n-2}\\
y_{n-2}\\
z_{n-2}
\end{matrix}
\right]+
\left[
\begin{matrix}
-x_{0}+y_{0}-z_{0}\\
-x_{0}+y_{0}-z_{0}\\
-x_{0}+y_{0}-z_{0}
\end{matrix}
\right]
\end{align*}
where in  the  last equality $(i)$ was used. This concludes the proof of Proposition \ref{proposition:const}.

\end{proof}

Noting that $C$ is invertible, the preceding recurrence relations allows us to define $x_{n},y_{n},z_{n}$ for negative values with the relation,
\begin{equation}\label{sist:recursao}
x_{n-3}=4x_{n-2}-4x_{n-1}+x_{n}
\end{equation}
and analogously for $y_{n}$ and $z_{n}$. We may, thus, define, for any $n \in \mathbb{Z}$,
\begin{equation*}
\left[
\begin{matrix}
x_{n}\\
y_{n}\\
z_{n}
\end{matrix}
\right]=
C^{n}
\left[
\begin{matrix}
a\\
b\\
c
\end{matrix}
\right]=
\left[
\begin{matrix}
2 &0 &-1\\
1 &0 &0\\
0 &-1 &2
\end{matrix}
\right]^{n}
\left[
\begin{matrix}
a\\
b\\
c
\end{matrix}
\right]
\end{equation*}

\bigbreak

We remark that the coefficients of $a, b, c$ in the leftmost color (as well as in the middle and in the rightmost colors), satisfy a recurrence relation similar to the one $x_{n}$ does.

\begin{definition}Set:
\begin{equation*}
\begin{matrix}
x_{n}=a_{n}a+b_{n}b+c_{n}c\\
y_{n}=a_{n}'a+b_{n}'b+c_{n}'c\\
z_{n}=a_{n}''a+b_{n}''b+c_{n}''c
\end{matrix}
\end{equation*}
\end{definition}

We have:
\begin{equation*}
\left[
\begin{matrix}
x_{n}\\
y_{n}\\
z_{n}
\end{matrix}
\right]=
C^{n}
\left[
\begin{matrix}
a\\
b\\
c
\end{matrix}
\right]=
\left[
\begin{matrix}
a_{n} &b_{n} &c_{n}\\
a_{n}' &b_{n}' &c_{n}'\\
a_{n}'' &b_{n}'' &c_{n}''
\end{matrix}
\right]
\left[
\begin{matrix}
a\\
b\\
c
\end{matrix}
\right]
\end{equation*}

\begin{corollary}
 Keeping the conditions of the proposition above:
\begin{enumerate}[(i)]
 \item $a_{n}=4a_{n-1}-4a_{n-2}+a_{n-3}$;\\
$b_{n}=4b_{n-1}-4b_{n-2}+b_{n-3}$;\\
$c_{n}=4c_{n-1}-4c_{n-2}+c_{n-3}$.

\item $a_{n}=3a_{n-1}-a_{n-2}-1$;\\
$b_{n}=3b_{n-1}-b_{n-2}+1$;\\
$c_{n}=3c_{n-1}-c_{n-2}-1$.
\end{enumerate}

\end{corollary}
\begin{proof}
 We will prove the proposition for the sequence $a_{n}$, since the other proofs are analogous. Consider the matrix $C^{n}$. Then, one realizes $a_{n}$ is the leftmost color after $\big( \sigma\sb{2}\sigma\sb{1}\sp{-1} \big)\sp{n}$, when $x_{0}=1, y_{0}=0$, $z_{0}=0$ i.e., $x\sb{n}=a\sb{n}$ when $x_{0}=1, y_{0}=0$, $z_{0}=0$.
 Applying Proposition \ref{proposition:const}, we obtain the desired relations.
\end{proof}

More generally, the entries of $C^{n}$ satisfy a recurrence relation similar to the one $x_{n}, y_{n}$,
and $z_{n}$ do, as $c_{ij}^{n}$, the entry $ij$ of $C^{n}$, can be interpreted as the colors of the $i-th$
strand with the top colors equal to zero except for the $j-th$ strand that takes 1 as initial color.

Before proceeding to the next result, let us look at the first few powers of $C$:

\begin{equation*}
C=
\left[
\begin{matrix}
2 &0 &-1\\
1 &0 &0\\
0 &-1 &2
\end{matrix}
\right],
\qquad
C^{2}=
\left[
\begin{matrix}
4 &1 &-4\\
2 &0 &-1\\
-1 &-2 &4
\end{matrix}
\right],
\qquad
C^{3}=
\left[
\begin{matrix}
9 &4 &-12\\
4 &1 &-4\\
-4 &-4 &9
\end{matrix}
\right]
\end{equation*}

\begin{corollary}\label{cor:matrix}
The powers of the matrix $C$ satisfy:

\begin{equation*}
C^{n}=
\left[
\begin{matrix}
a_{n} &b_{n} &-b_{n+1}\\
a_{n-1} &b_{n-1} &-b_{n}\\
-b_{n} &-a_{n-1} &a_{n}
\end{matrix}
\right], n\in \mathbb{Z}
\end{equation*}

\end{corollary}

\begin{proof}
We already saw that the entries of the matrices $C^{n}$ satisfy a specific recurrence relation, which is similar to the ones $a_{n}$, and $b_{n}$ do. So, it is enough to verify that the first 3 powers of $C$ satisfy the expression above, and to use induction to establish the result. We leave the details to the reader.
\end{proof}


We recall we want to establish a formula which yields the number of $r$-colorings of $THK(3, n)$. In order to do that, we solve the following system of linear equations over $\mathbb{Z}_{r}$:
\begin{equation}\label{sistema}
C\sp{n}\left[
\begin{matrix}
x_{0}\\
y_{0}\\
z_{0}
\end{matrix}
\right]=
\left[
\begin{matrix}
x_{0}\\
y_{0}\\
z_{0}
\end{matrix}
\right]
\end{equation}

which, upon rewriting and applying Corollary \ref{cor:matrix} yields the following system of linear homogeneous equations over $\mathbb{Z}\sb{r}$:

\begin{equation}\label{sistema2}
\left[
\begin{matrix}
a_{n}-1 &b_{n} &-b_{n+1}\\
a_{n-1} &b_{n-1}-1 &-b_{n}\\
-b_{n} &-a_{n-1} &a_{n}-1
\end{matrix}
\right]
\left[
\begin{matrix}
x_{0}\\
y_{0}\\
z_{0}
\end{matrix}\right]
=\left[
\begin{matrix}
0\\
0\\
0
\end{matrix}
\right]
\end{equation}

This coefficient matrix in (\ref{sistema2}) can be further simplified as we will now show.

\begin{corollary}\label{corollary:index}
For any positive integer $n$:
\begin{align*}
 a_{n}-a_{n-1}-b_{n}&=1\\
b_{n}-b_{n-1}-a_{n-1}&=-1
\end{align*}
\end{corollary}

\begin{proof}
 We remark that $a_{n},a_{n-1}$ and $-b_{n}$ are the bottom colors, from left to right, respectively, of the
 braid $(\sigma_{2} \sigma_{1}^{-1})^{n}$ when $x_{0}=1,y_{0}=0$ and $z_{0}=0$. Therefore, by Proposition \ref{proposition:const},
 we have:

\begin{equation*}
a_{n}-a_{n-1}+(-b_{n})=1-0+0=1
\end{equation*}

Also, $b_{n}, b_{n-1}$ and $-a_{n-1}$ are the bottom colors, from left to right, respectively, of the strings of the braid $(\sigma_{2} \sigma_{1}^{-1})^{m}$ when $x_{0}=0,y_{0}=1$ and $z_{0}=0$. Again by Proposition \ref{proposition:const}:

\begin{equation*}
b_{n}-b_{n-1}+(-a_{n-1})=0-1+0=-1
\end{equation*}
\end{proof}

Now, by Corollary \ref{corollary:index} we may conclude that:

\begin{equation*}
\left[
\begin{matrix}
1 &-1 &1
\end{matrix}
\right]
\left[
\begin{matrix}
a_{n} &b_{n} &-b_{n+1}\\
a_{n-1} &b_{n-1} &-b_{n}\\
-b_{n} &-a_{n-1} &a_{n}
\end{matrix}
\right]
-\left[
\begin{matrix}
1 &-1 &1
\end{matrix}
\right]=\left[
\begin{matrix}
0 &0 &0
\end{matrix}
\right]
\end{equation*}

In this way, by adding the first line and subtracting the second line to the third line in $C^{n}-I_{3}$,
we obtain:

\begin{equation*}
\left[
\begin{matrix}
1 &0 &0\\
0 &1 &0\\
1 &-1 &1
\end{matrix}
\right]
\left[
\begin{matrix}
a_{n}-1 &b_{n} &-b_{n+1}\\
a_{n-1} &b_{n-1}-1 &-b_{n}\\
-b_{n} &-a_{n-1} &a_{n}-1
\end{matrix}
\right]=
\left[
\begin{matrix}
a_{n}-1 &b_{n} &-b_{n+1}\\
a_{n-1} &b_{n-1}-1 &-b_{n}\\
0 &0 &0
\end{matrix}
\right]
\end{equation*}


We now define the sequences $u_{n}$ and $v_{n}$ which will relate to $a_{n}$ and $b_{n}$.

\begin{definition}
Let $u_{n}$ and $v_{n}$ be sequences defined recursively as follows:

\begin{equation*}
\begin{cases}
u_{-3}=-1,\\
u_{-2}=-1,\\
u_{-1}=0,\\
u_{0}=1,\\
u_{n}=3u_{n-2}-u_{n-4}, &\text{if } n \geq 1
\end{cases}
\end{equation*}
and,
\begin{equation*}
\begin{cases}
v_{-3}=7,\\
v_{-2}=2,\\
v_{-1}=3,\\
v_{0}=1,\\
v_{n}=3v_{n-2}-v_{n-4}, &\text{if } n \geq 1
\end{cases}
\end{equation*}

\end{definition}

\begin{proposition}\label{recoprod}
For $n \in \mathbb{N}$ we have:
\begin{alignat*}{2}
a_{n}&=u_{n}v_{n} & b_{n}&=u_{n-2}u_{n-1}\\
(a_{n}-1)&=u_{n-1}v_{n+1} \qquad & (b_{n}-1)&=u_{n}u_{n-3}
\end{alignat*}

\end{proposition}

\begin{proof}
First we need to verify the cases $n=0,1,2,3,4$. We leave this task to the reader. We will prove by induction that
 $a_{n}=u_{n}v_{n}$ and $(a_{n}-1)=u_{n-1}v_{n+1}$. For $n\geq 0$ we assume the validity of the hypothesis for
 $n,n+1,n+2,n+3$ and $n+4$ and conclude that it is also valid for $n+5$. Proposition \ref{proposition:const} will be used throughout the following calculations.

\begin{align*}
u_{n+5}v_{n+5}&=(3u_{n+3}-u_{n+1})(3v_{n+3}-v_{n+1})=\\
&=9u_{n+3}v_{n+3}-3(3u_{n+1}-u_{n-1})v_{n+1}-3u_{n+1}v_{n+3}+u_{n+1}v_{n+1}=\\
&=9a_{n+3}-9a_{n+1}+3(a_{n}-1)-3(a_{n+2}-1)+a_{n+1}=\\
&=9a_{n+3}-15a_{n+2}+4a_{n+1}+3(4a_{n+2}-4_{n+1}+a_{n})=\\
&=-4a_{n+3}+a_{n+2}+4(4a_{n+3}-4a_{n+2}+a_{n+1})=\\
&=4a_{n+4}-4a_{n+3}+a_{n+2}=a_{n+5}
\end{align*}

\begin{align*}
u_{n+4}v_{n+6}&=(3u_{n+2}-u_{n})(3v_{n+4}-v_{n+2})=\\
&=9u_{n+2}v_{n+4}-3u_{n}(3v_{n+2}-v_{n})-3u_{n+2}v_{n+2}+u_{n}v_{n+2}=\\
&=9(a_{n+3}-1)-9(a_{n+1}-1)+3a_{n}-3a_{n+2}+(a_{n+1}-1)=\\
&=9a_{n+3}-15a_{n+2}+4a_{n+1}-1+3(4a_{n+2}-4a_{n+1}+a_{n})=\\
&=-4a_{n+3}+a_{n+2}-1+4(4a_{n+3}-4a_{n+2}+a_{n+1})=\\
&=4a_{n+4}-4a_{n+3}+a_{n+2}-1=a_{n+5}-1
\end{align*}

This ends the first part of the proof. Now, let us prove by induction that $b_{n}=u_{n-2}u_{n-1}$ and
$(b_{n}-1)=u_{n}u_{n-3}$. We assume again for $n\geq 0$ the validity of the hypothesis for $n,n+1,n+2,n+3$ and
$n+4$, concluding that it is valid for $n+5$.

\begin{align*}
u_{n+3}u_{n+4}&=u_{n+3}(3u_{n+2}-u_{n})=3b_{n+4}-(b_{n+3}-1)=\\
&=3b_{n+4}-b_{n+3}+1=b_{n+5}
\end{align*}

\begin{align*}
u_{n+5}u_{n+2}&=(3u_{n+3}-u_{n+1})u_{n+2}=3b_{n+4}-b_{n+3}=\\
&=(3b_{n+4}-b_{n+3}+1)-1=b_{n+5}-1
\end{align*}

This ends the proof.

\end{proof}

The linear homogeneous system of equations over $\mathbb{Z}\sb{r}$ is now equivalent to:
\begin{equation}\label{eqn:lhse}
u_{n-1}
\left[
\begin{matrix}
v_{n+1} &u_{n-2} &-u_{n}\\
v_{n-1} &u_{n-4} &-u_{n-2}\\
0 &0 &0
\end{matrix}
\right]\left[
\begin{matrix}
x_{0}\\
y_{0}\\
z_{0}
\end{matrix}
\right]=\left[
\begin{matrix}
0\\
0\\
0
\end{matrix}
\right]
\end{equation}

By performing elementary operations on the lines of the matrix we obtain an even simpler coefficient matrix:
\begin{multline*}
\left[
\begin{matrix}
v_{k+1} &u_{k-2} &-u_{k}\\
v_{k-1} &u_{k-4} &-u_{k-2}\\
0 &0 &0
\end{matrix}
\right]\longrightarrow
\left[
\begin{matrix}
v_{k+1}-2v_{k-1} &u_{k-2}-2u_{k-4} &-u_{k}+2u_{k-2}\\
v_{k-1} &u_{k-4} &-u_{k-2}\\
0 &0 &0
\end{matrix}
\right]\longrightarrow\\
\\
\longrightarrow\left[
\begin{matrix}
v_{k+1}-2v_{k-1} &u_{k-2}-2u_{k-4} &-u_{k}+2u_{k-2}\\
-v_{k+1}+3v_{k-1} &-u_{k-2}+3u_{k-4} &u_{k}-3u_{k-2}\\
0 &0 &0
\end{matrix}
\right]
=\\
\\
=\left[
\begin{matrix}
v_{k-1}-v_{k-3} &u_{k-4}-u_{k-6} &-u_{k-2}+u_{k-4}\\
v_{k-3} &u_{k-6} &-u_{k-4}\\
0 &0 &0
\end{matrix}
\right]
\longrightarrow\left[
\begin{matrix}
v_{k-1} &u_{k-4} &-u_{k-2}\\
v_{k-3} &u_{k-6} &-u_{k-4}\\
0 &0 &0
\end{matrix}
\right]
\end{multline*}

Therefore, for $n$ odd, \ref{eqn:lhse} simplifies to:

\begin{equation}\label{mat:odd}
u_{n-1}\left[
\begin{matrix}
v_{2} &u_{-1} &-u_{1}\\
v_{0} &u_{-3} &-u_{-1}\\
0 &0 &0
\end{matrix}
\right]=u_{n-1}\left[
\begin{matrix}
1 &0 &-1\\
1 &-1 &0\\
0 &0 &0
\end{matrix}
\right]\longrightarrow u_{n-1}\left[
\begin{matrix}
1 &0 &-1\\
0 &-1 &1\\
0 &0 &0
\end{matrix}
\right]
\end{equation}

whereas, for $n$ even, \ref{eqn:lhse} simplifies to:
\begin{equation}\label{mat:even}
u_{n-1}\left[
\begin{matrix}
v_{3} &u_{0} &-u_{2}\\
v_{1} &u_{-2} &-u_{0}\\
0 &0 &0
\end{matrix}
\right]=u_{n-1}\left[
\begin{matrix}
3 &1 &-4\\
2 &-1 &-1\\
0 &0 &0
\end{matrix}
\right]\longrightarrow u_{n-1}\left[
\begin{matrix}
1 &2 &-3\\
0 &-5 &5\\
0 &0 &0
\end{matrix}
\right]
\end{equation}

We can now establish the formulas for the number of the colorings in terms of $r$ and $u\sb{n-1}$ in Theorem \ref{thm:numbercols}. We state that part here again as Proposition \ref{prop:numbercols} below for the reader's convenience. We recall we let $(a, b)$ stand for the greatest common divisor of the positive integers $a$ and $b$.

\begin{proposition}\label{prop:numbercols} Given positive integers $n$ and $r>1$,
\[
\# col\sb{r} THK(3, n) =
\begin{cases}
(u\sb{n-1}, r)\sp{2}r & \text{ if $n$ is odd}\\
(5u\sb{n-1}, r)(u\sb{n-1}, r)r & \text{ if $n$ is even}
\end{cases}
\]
\end{proposition}
\begin{proof}Let us suppose that $n$ is odd. The color input $a,b,c$ in $\big( \sigma\sb{2}\sigma\sb{1}\sp{-1} \big)\sp{n}$ induces an $r$-coloring in $THK(3, n)$ if and only if:
\begin{equation*}
u_{n-1}\left[
\begin{matrix}
1 &0 &-1\\
0 &-1 &1\\
0 &0 &0
\end{matrix}
\right]
\left[
\begin{matrix}
a\\
b\\
c
\end{matrix}
\right]=_{r}
\left[
\begin{matrix}
0\\
0\\
0
\end{matrix}
\right]
\end{equation*}

We now count the solutions of this linear homogeneous system of equations over $\mathbb{Z}_{r}$ (see \cite{kl}, Claims 2.1 and 2.2). By inspection of the third equation, we realize $c$ can take on any value from $\mathbb{Z}_{r}$ i.e., $r$ possibilities for $c$. The second equation, $u_{n-1}(-b+c)=_{r} 0$, yields $b-c=k\frac{r}{(u_{n-1},r)}$ for $k=0,1,...,(u_{n-1},r)-1$. Note that for each positive integer $a$,
\[
ax=\sb{r}0 \, \Leftrightarrow \, (a, r)x=\sb{r}0
\]
Moreover, from the first equation, $a=c+k\frac{r}{(u_{n-1},r)}$. So, both $b$ and $a$ have $(u_{n-1},r)$ possibilities, compatible with each of the $r$ possibilities for $c$. We conclude that the indicated system of equations has $r(u_{n-1},r)^{2}$ solutions over $\mathbb{Z}_{r}$.

Let us now suppose that $n$ is even. The color input $a,b,c$ in $\big( \sigma\sb{2}\sigma\sb{1}\sp{-1} \big)\sp{n}$ induces an $r$-coloring in $THK(3, n)$ if and only if:
\begin{equation*}
u_{n-1}\left[
\begin{matrix}
1 &2 &-3\\
0 &-5 &5\\
0 &0 &0
\end{matrix}
\right]
\left[
\begin{matrix}
a\\
b\\
c
\end{matrix}
\right]
\equiv_{r}
\left[
\begin{matrix}
0\\
0\\
0
\end{matrix}
\right]
\end{equation*}

Arguing as above, we conclude that the number of solutions over $\mathbb{Z}_{r}$ of the indicated system of equations is $r(5u_{n-1},r)(u_{n-1},r)$, for even $n$. This concludes the proof.
\end{proof}

\bigbreak

\begin{proposition}
For $n \in \mathbb{Z}$,

\begin{equation}\label{bu}
u\sb{n} = \frac{1}{\sqrt{5}}\left(\left(\frac{1+\sqrt{5}}{2}\right)^{n+2}-\left(\frac{-1+\sqrt{5}}{2}\right)^{n}-
\left(\frac{1-\sqrt{5}}{2}\right)^{n+2}+\left(\frac{-1-\sqrt{5}}{2}\right)^{n}\right)
\end{equation}
\end{proposition}
\begin{proof}
The $u\sb{n}$ is a solution of the linear difference equation with the ``initial values'':
\[
-u\sb{n}+3u\sb{n-2}-u\sb{n-4}=0 \qquad \qquad u\sb{-3}=-1 \qquad u\sb{-2}=-1  \qquad u\sb{-1}=0  \qquad u\sb{0}=1
\]
The general method of solving this type of initial value problem can be found in \cite{Henrici}, for instance.

The characteristic equation here is
\[
x^{4}-3x^{2}+1=0
\]
so:

\begin{equation*}
 x^{2}=\left(\frac{3 \pm \sqrt{5}}{2}\right)=\left(\frac{1+5\pm 2\sqrt{5}}{4}\right)
 \Leftrightarrow x=\pm\left(\frac{1\pm\sqrt{5}}{2}\right)
\end{equation*}

Then, the general solution of $-u\sb{n}+3u\sb{n-2}-u\sb{n-4}=0$ is:
\[
u\sb{n}=c\sb{1}\left(\frac{1+\sqrt{5}}{2}\right)^{n}+c\sb{2}\left(\frac{-1+\sqrt{5}}{2}\right)^{n}+
c\sb{3}\left(\frac{1-\sqrt{5}}{2}\right)^{n}+c\sb{4}\left(-\frac{1-\sqrt{5}}{2}\right)^{n}
\]

By using the ``initial values'' $u\sb{-3}=-1, u\sb{-2}=-1, u\sb{-1}=0, u\sb{0}=1$, we obtain a system of four linear equations on the $c\sb{i}$'s.  Solving this system of equations we obtain the $c\sb{i}$'s and finally:
\[
u\sb{n}= \frac{1}{\sqrt{5}}\left(\left(\frac{1+\sqrt{5}}{2}\right)^{n+2}-\left(\frac{-1+\sqrt{5}}{2}\right)^{n}-
\left(\frac{1-\sqrt{5}}{2}\right)^{n+2}+\left(\frac{-1-\sqrt{5}}{2}\right)^{n}\right)
\]
\end{proof}

\bigbreak

This concludes the proof of Theorem \ref{thm:numbercols} and of Corollary \ref{cor:nontriv}.

The formulas for the determinants in Corollary \ref{cor:det} follow from the observation of matrices (\ref{mat:odd}) for the odd $n$ case and (\ref{mat:even}) for the even $n$ case; the fact that these determinants are always greater than zero follows from proving by induction that
\[
u\sb{n} > 0 \quad \text{ and } \quad u\sb{n}-u\sb{n-2} > 0 \qquad \text{ for all } n>2
\]

\bigbreak

\subsection{Proof of Theorem \ref{thm:minexact}}\label{subsect:thmminexact}

\noindent

\subsubsection{Preliminaries}\label{subsubsect:prelim}

\noindent

We first establish Propositions \ref{prop:pcolrcol}, \ref{prop:stackcol}, \ref{prop:nonsplit} for they will be useful in the sequel.

\begin{proposition}\label{prop:pcolrcol} Suppose a knot admits a non-trivial $s$-coloring and $s\mid r$. Then this knot also admits a non-trivial $r$-coloring.
\end{proposition}\begin{proof} Since $s\mid r$, then the set
\[
\bigg\{\,  0, \frac{r}{s}, 2\frac{r}{s}, \dots , (s-1)\frac{r}{s}  \, \bigg\} \quad \text{ mod } r
\]
is closed with respect to the $a\ast b:= 2b-a $ mod $r$ operation. Moreover, the mapping
\begin{align*}
f \, : \quad &\mathbb{Z}\sb{s} \longrightarrow \mathbb{Z}\sb{r}\\
&\, i \longrightarrow i\frac{r}{s}
\end{align*}
is injective and preserves the $\ast $ operation. In this way, if $(i\sb{1}, i\sb{2}, \dots , i\sb{N})$ is the sequence of colors mod $s$ one should assign to the sequence of arcs in a knot diagram of the knot under study to obtain a non-trivial $s$-coloring, then the sequence $(i\sb{1}\frac{r}{s}, i\sb{2}\frac{r}{s}, \dots , i\sb{N}\frac{r}{s})$ assigned to the same sequence of arcs, represents a non-trivial $r$-coloring of the same knot.
\end{proof}

\bigbreak

\begin{proposition}\label{prop:stackcol} Consider the positive integers $c, k, n$, and integer $r>1$. Suppose the standard diagram of $THK(3, n)$ admits a non-trivial $r$-coloring with $c$ colors. Then the standard diagram of $THK(3, kn)$ also admits a non-trivial $r$-coloring with $c$ colors.
\end{proposition}\begin{proof} If the standard diagram of $THK(3, n)$ admits a non-trivial $r$-coloring with $c$ colors, we consider this coloring in $\big( \sigma\sb{2}\sigma\sb{1}\sp{-1}  \big)\sp{n}$ i.e., before braid closure. In this way, the sequence of colors on the top arcs (from left to right) equals the sequence of colors on the bottom arcs (from left to right). We then juxtapose $k$ copies of this colored $\big( \sigma\sb{2}\sigma\sb{1}\sp{-1}  \big)\sp{n}$. Upon taking its closure we obtain a non-trivial $r$-coloring of $THK(3, kn)$.
\end{proof}

\bigbreak

\begin{definition} A knot is said {\rm split} if there exist two disjoint neighborhoods in $3$-space, say $N\sb{1}$ and $N\sb{2}$, such that it is deformable into a knot, say $K$, such that
\[
K \subset N\sb{1} \cup N\sb{2} \qquad K\cap N\sb{1} \neq \emptyset \qquad K\cap N\sb{2} \neq \emptyset
\]
Otherwise, the knot is said {\rm non-split}.
\end{definition}

\begin{proposition}\label{prop:nonsplit} For any positive integer $n$, $THK(3, n)$ is non-split.
\end{proposition}
\begin{proof} Fix a positive integer $n$. If $THK(3, n)$ were split, then for any integer $r>1$, there would be at least $r\sp{2}$ $r$-colorings. Each of these $r\sp{2}$ $r$-colorings stands for the $r$-coloring which results from trivially $r$-coloring each of the two neighborhoods the knot splits into. But, setting $r$ equal to a prime larger than both $u\sb{n-1}$ and $5$, yields
\[
\# col\sb{r} THK(3, n) \quad = \quad
\begin{cases}
(u\sb{n-1}, r)\sp{2}r & \text{ if $n$ is odd}\\
(5u\sb{n-1}, r)(u\sb{n-1}, r)r & \text{ if $n$ is even}
\end{cases}
\quad = \quad r < r\sp{2}
\]
Hence, $THK(3, n)$ is non-split.
\end{proof}

\subsubsection{Further analysis of the $u\sb{n}$ sequence and the $\psi $ mapping: towards the proof of Theorem \ref{thm:minexact}}\label{subsubsect:psi}

\noindent

In order to prove the ``if'' parts in Theorem \ref{thm:minexact}, and to prove Theorem \ref{thm:minestimates}, we analyze further the sequence $u\sb{n}$.

\bigbreak

\begin{proposition}\label{proposition:sominhas}
For $n \in \mathbb{Z}$ we have:
\begin{align*}
u_{2n}&=u_{2n+1}+u_{2n-1}\\
u_{2n+1}&=\frac{u_{2n+2}+u_{2n}}{5}
\end{align*}
\end{proposition}

\begin{proof}
We only prove the inductive step:
\begin{align*}
u\sb{2n+3}+u\sb{2n+1}&=3u\sb{2n+1}-u\sb{2n-1}+u\sb{2n+1}=5u\sb{2n+1}-u\sb{2n+1}-u\sb{2n-1}=\\
&=u\sb{2n+2}+(u\sb{2n}-u\sb{2n+1}-u\sb{2n-1}) = u\sb{2n+2}+0
\end{align*}
\begin{align*}
u\sb{2n+4}+u\sb{2n+2}&=3u\sb{2n+2}-u\sb{2n}+u\sb{2n+2}=5u\sb{2n+2}-u\sb{2n+2}-u\sb{2n}=\\
&=5(u\sb{2n+3}+u\sb{2n+1})-u\sb{2n+2}-u\sb{2n} =5u\sb{2n+3}+5(u\sb{2n+1}-u\sb{2n+2}-u\sb{2n})= 5u\sb{2n+3}+0
\end{align*}
\end{proof}

\begin{proposition}\label{proposition:produ}
\begin{enumerate}[(i)]
\item If $m$ is even, or $n$ is odd, then:
\begin{equation}\label{eqprod}
u_{m+n}=u_{m+1}u_{n}-u_{m-1}u_{n-2}
\end{equation}

\item If $m$ is even and $n$ is odd, then:

\begin{equation*}
u_{m+n}=u_{m}u_{n}-u_{m-1}u_{n-1}
\end{equation*}
\end{enumerate}

\end{proposition}

\begin{proof}
To prove $(i)$ let us first suppose $m$ is an arbitrary integer, we have:

\begin{align*}
&u_{m+1}=u_{m+1}\times 1-u_{m-1}\times 0=u_{m+1}u_{1}-u_{m-1}u_{-1}\\
&u_{m+3}=u_{m+1}\times 3-u_{m-1}\times 1=u_{m+1}u_{3}-u_{m-1}u_{1}
\end{align*}

Also, if $m$ is even, by Proposition \ref{proposition:sominhas}, we have:

\begin{align*}
&u_{m}=u_{m+1}\times 1-u_{m-1}\times(-1)=u_{m+1}u_{0}-u_{m-1}u_{-2}\\
&u_{m+2}=u_{m+3}+u_{m+1}=u_{m+1}\times 4-u_{m-1}\times 1=u_{m+1}u_{2}-u_{m-1}u_{0}
\end{align*}

Now, fixed $m$, we will do induction on $n$. Let us suppose that the equation (\ref{eqprod}) is true for a given $m$ and $n=k,k+2$, then:

\begin{align*}
u_{m+(k+4)}&=3u_{m+(k+2)}-u_{m+k}=u_{m+1}(3u_{k+2}-u_{k})-u_{m-1}(3u_{k}-u_{k-2})=\\
&=u_{m+1}u_{k+4}-u_{m-1}u_{k+2}\\
u_{m+(k-2)}&=3u_{m+k}-u_{m+(k+2)}=u_{m+1}(3u_{k}-u_{k+2})-u_{m-1}(3u_{k-2}-u_{k})=\\
&=u_{m+1}u_{k-2}-u_{m-1}u_{k-4}
\end{align*}

So, equation (\ref{eqprod}) is also verified for $n=k-2,k+4$. From here we may conclude that when $m$ is even equation (\ref{eqprod}) is true for every $n$, and when $m$ is odd equation (\ref{eqprod}) is true for $n$ odd.

To prove $(ii)$ we use $(i)$, and Proposition \ref{proposition:sominhas}:

\begin{align*}
u_{m+n}&=u_{m+1}u_{n}-u_{m-1}u_{n-2}=(u_{m}-u_{m-1})u_{n}-u_{m-1}u_{n-2}\\
&=u_{m}u_{n}-u_{m-1}(u_{n}+u_{n-2})=u_{m}u_{n}-u_{m-1}u_{n-1}
\end{align*}
\end{proof}

\bigbreak

We repeat here the definition of the $\psi$ mapping for the reader's convenience.

\begin{definition}\label{def:psi2} For any integer $r>1$ set
\[
\psi (r) := \min \{ \, q\in \mathbb{Z}\sp{+}\, \big| \, r\mid u\sb{q-1}  \,   \}
\]
\end{definition}

\bigbreak

\begin{proposition}
$\psi$ is well-defined.
\end{proposition}

\begin{proof}
We will prove that, for any integer $r>1$, $\{ q \in \mathbb{Z}^{+}\, \big|\,  r\mid u_{q-1}\}\neq \emptyset$. Since it is a subset of $\mathbb{Z}\sp{+}$, then it has a minimum.

Fix an integer $r>1$ and set $U_{n}=(u_{n},u_{n+1},u_{n+2},u_{n+3})$ a sequence in $\mathbb{Z}_{r}^{4}$ formed by consecutive terms of the $u$ sequence read mod $r$. Let us consider the first $r^{4}+1$  $U_{i}$'s, $i=0, \dots ,r^{4}$. Since $\mathbb{Z}_{r}^{4}$ has $r^{4}$ elements, by the Pigeonhole Principle, there are integers $i,j$, such that $U_{i}=U_{j}$, $0\leq i< j \leq r^{4}$. Hence, as $u_{n}$ can be defined recursively by the previous four terms of the same sequence, we conclude that $\{u_{n} \}_{n\in \mathbb{Z}\sp{+}}$ (mod $r$) is periodic with period $j-i$ (or less). So, $u_{j-i-1}\equiv_{r} u_{-1}\equiv_{r} 0$ and $\{ q \in \mathbb{Z}^{+}\mid r\mid u_{q-1}\}\ni (j-i)$ is non-empty, ending the proof.
\end{proof}

\begin{proposition}\label{proposition:divide}
$r \mid u_{m-1}$ if and only if $\psi(r)\mid m$.
\end{proposition}

\begin{proof} Upon extending $u\sb{n}$ to the negative integers using the recurrence relation $u\sb{n-4}=3u\sb{n-2}-u\sb{n}$, we note that $u\sb{n}=-u\sb{-n-2}$.

We start by verifying that $u_{n}=-u_{-n-2}$. This is done using the fact that $u_{1}=-u_{-3}=1$, $u_{0}=-u_{-2}=1$ and $u_{-1}=-u_{-1}=0$. After this, we leave it as an exercise to use induction on $n$ to obtain the desired conclusion.

Now, let us suppose that $u_{\psi(n)-1} \mid u_{m-1}$, to conclude that $u_{\psi(n)-1} \mid u_{m \pm \psi(n)-1}$.

If $\psi(n)$ is even or $(m-1)$ is odd, by Proposition \ref{proposition:produ} we have:
\begin{align*}
u_{\pm \psi(n)+(m-1)}&=u_{(\pm \psi(n)+1)}u_{m-1}-u_{\pm \psi(n)-1}u_{m-3}=\\
&=u_{(\pm \psi(n)+1)}u_{m-1}\mp u_{\psi(n)-1}u_{m-3}
\end{align*}

Otherwise, if $\psi(n)$ is odd and $(m-1)$ is even, by Proposition \ref{proposition:produ} we have:
\begin{align*}
u_{(m-1)\pm \psi(n)}&=u_{(m-1)}u_{\pm \psi(n)}-u_{m-2}u_{\pm \psi(n)-1}=\\
&=u_{(m-1)}u_{\pm \psi(n)}\mp u_{m-2}u_{\psi(n)-1}
\end{align*}

Either way, we conclude that $u_{\psi(n)-1} \mid u_{m \pm \psi(n)-1}$.

From what we have seen above, we obtain, by induction, that $n\mid u_{\psi(n)-1} \mid u_{k\psi(n)-1}$,
for $k\in \mathbb{Z}$. Also, if $r$ is the remainder of the division of $m$ by $\psi(n)$, we conclude that $n \mid u_{\psi(n)-1}\mid u_{r-1}$. Hence, by definition of $\psi(n)$, we get $r=0$, concluding the proof.
\end{proof}

\bigbreak

\subsubsection{The proof of Theorem \ref{thm:minexact}}\label{subsubsect:thmminexact}

\noindent

The proof of Theorem \ref{thm:minexact} now follows easily thanks to the result in \cite{lm} which we reproduce here for the reader's convenience. We remark that, given two positive integers $a, b$, we let $\langle a, b\rangle$ stand for $1$ if $a, b$ are relatively prime, otherwise, we let  $\langle a, b\rangle$ stand for their least common prime factor.

\begin{theorem}\label{thm:saito}Let $r$ be an integer greater than $1$.

Let $L$ be a non-split link.

\begin{enumerate}
\item  \quad $\langle r, \det L \rangle = 2$  \quad if, and only if,  \quad $mincol\sb{r}L = 2$
\item   \quad $\langle r, \det L \rangle = 3$  \quad if, and only if,  \quad $mincol\sb{r}L = 3$

\bigbreak

Furthermore, if $\det L \neq 0$ then:

\bigbreak

\item   \quad $\langle r, \det L \rangle \in \{ 5, 7 \}$  \quad if, and only if,  \quad $mincol\sb{r}L = 4$
\item    \quad $\langle r, \det L \rangle > 7$, \quad if, and only if,  \quad $mincol\sb{r}L \geq 5$
\end{enumerate}
\end{theorem}

\bigbreak
\begin{proof}(of Theorem \ref{thm:minexact})

\begin{enumerate}
\item Knowing that $\psi (2) = 3$,
\begin{align*}
&2\mid r \text{ and } 3\mid n \quad \Leftrightarrow \quad 2\mid r \text{ and } \psi (2)\mid n  \quad \Leftrightarrow  \quad 2\mid r \text{ and } 2\mid u\sb{n-1}  \quad \Leftrightarrow  \quad \\
&2\mid r \text{ and } 2\mid \det THK (3, n)  \quad \Leftrightarrow  \quad \langle r, \det THK (3, n) \rangle = 2 \quad \Leftrightarrow  \quad \\
&mincol\sb{r}THK (3, n) = 2
\end{align*}
\item Knowing that $2\nmid r \text{ or } 3\nmid n$ and $\psi (3)=4$,
\begin{align*}
&3\mid r \text{ and } 4\mid n \quad \Leftrightarrow \quad 3\mid r \text{ and } \psi (3)\mid n  \quad \Leftrightarrow  \quad 3\mid r \text{ and } 3\mid u\sb{n-1}  \quad \Leftrightarrow  \quad \\
&3\mid r \text{ and } 3\mid \det THK (3, n)  \quad \Leftrightarrow  \quad \langle r, \det THK (3, n) \rangle = 3 \quad \Leftrightarrow  \quad \\
&mincol\sb{r}THK (3, n) = 3
\end{align*}
\item The argument for this instance mimics the ones used in the preceding two instances. We leave the details to the reader.
\item The right-hand side of Figure \ref{fig:thmminexact2} shows an $11$-coloring of $THK(3, n)$ with $5$ colors. Hence, Proposition \ref{prop:pcolrcol}, Proposition \ref{prop:stackcol} along with the preceding instances, imply that if $11\mid r$ and $5\mid n$ (and neither $r$ nor $n$ comply with the preceding instances) then $mincol\sb{r}THK(3, n) = 5$.
\end{enumerate}
\end{proof}

\bigbreak

\subsection{Proof of Theorem \ref{thm:minestimates}}\label{subsect:thmminestimates}

\noindent

\subsubsection{Preliminaries}\label{subsubsect:prelim2}

\noindent

\begin{figure}[!ht]
    \psfrag{0}{\huge $0$}
    \psfrag{1}{\huge $1$}
    \psfrag{2}{\huge $2$}
    \psfrag{3}{\huge $3$}
    \psfrag{4}{\huge $4$}
    \psfrag{7}{\huge $7$}
    \psfrag{2-}{\huge A $2$-coloring of $THK(3, 3)$ with $2$ colors}
    \psfrag{5-}{\huge A $5$-coloring of $THK(3, 2)$ with $4$ colors}
    \centerline{\scalebox{.50}{\includegraphics{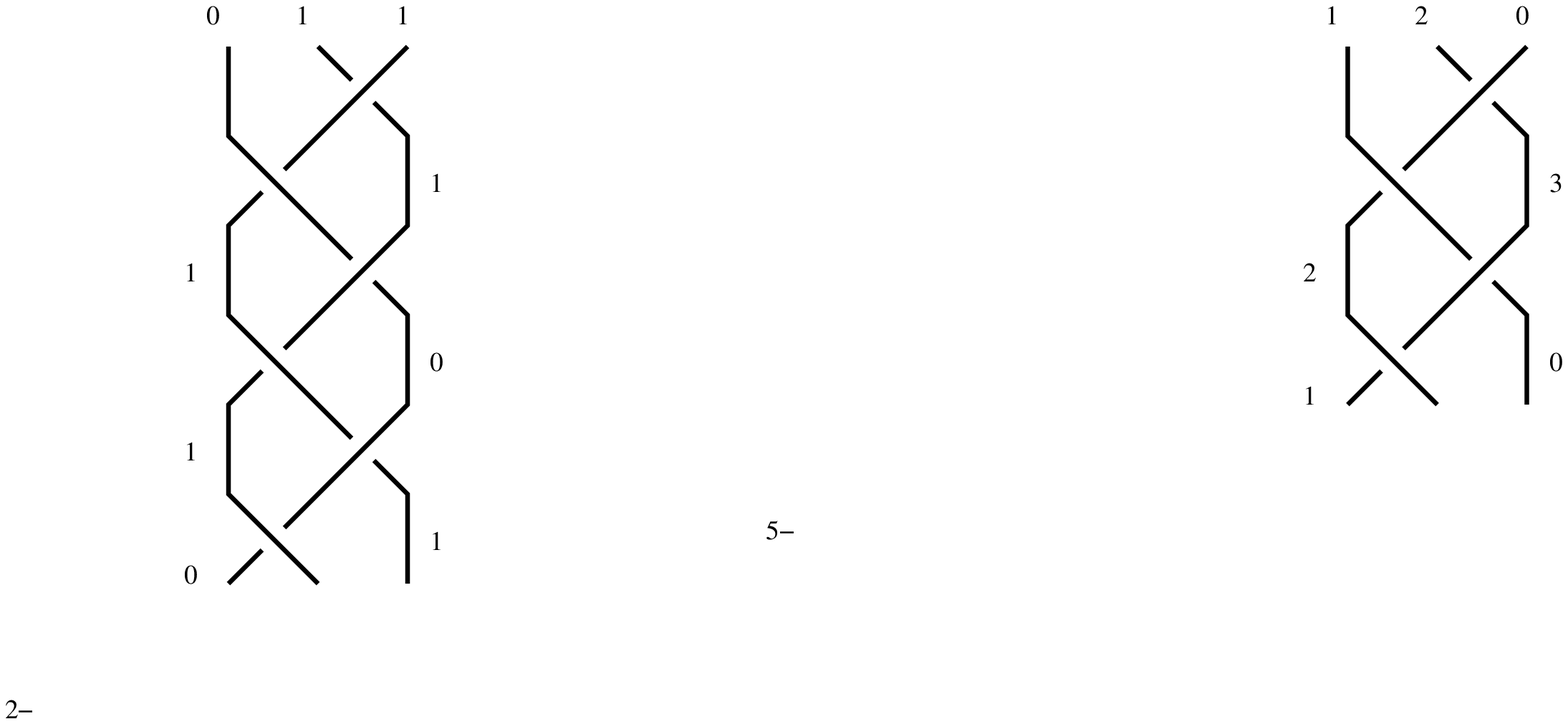}}}
    \caption{Illustrative examples of non-trivial colorings of $THK(3, n)$'s}\label{fig:thmminexact1}
\end{figure}

\begin{figure}[!ht]
    \psfrag{0}{\huge $0$}
    \psfrag{1}{\huge $1$}
    \psfrag{2}{\huge $2$}
    \psfrag{3}{\huge $3$}
    \psfrag{4}{\huge $4$}
    \psfrag{7}{\huge $7$}
    \psfrag{3-}{\huge A $3$-coloring of $THK(3, 4)$ with $3$ colors}
    \psfrag{11-}{\huge An $11$-coloring of $THK(3, 5)$ with $5$ colors}
    \centerline{\scalebox{.50}{\includegraphics{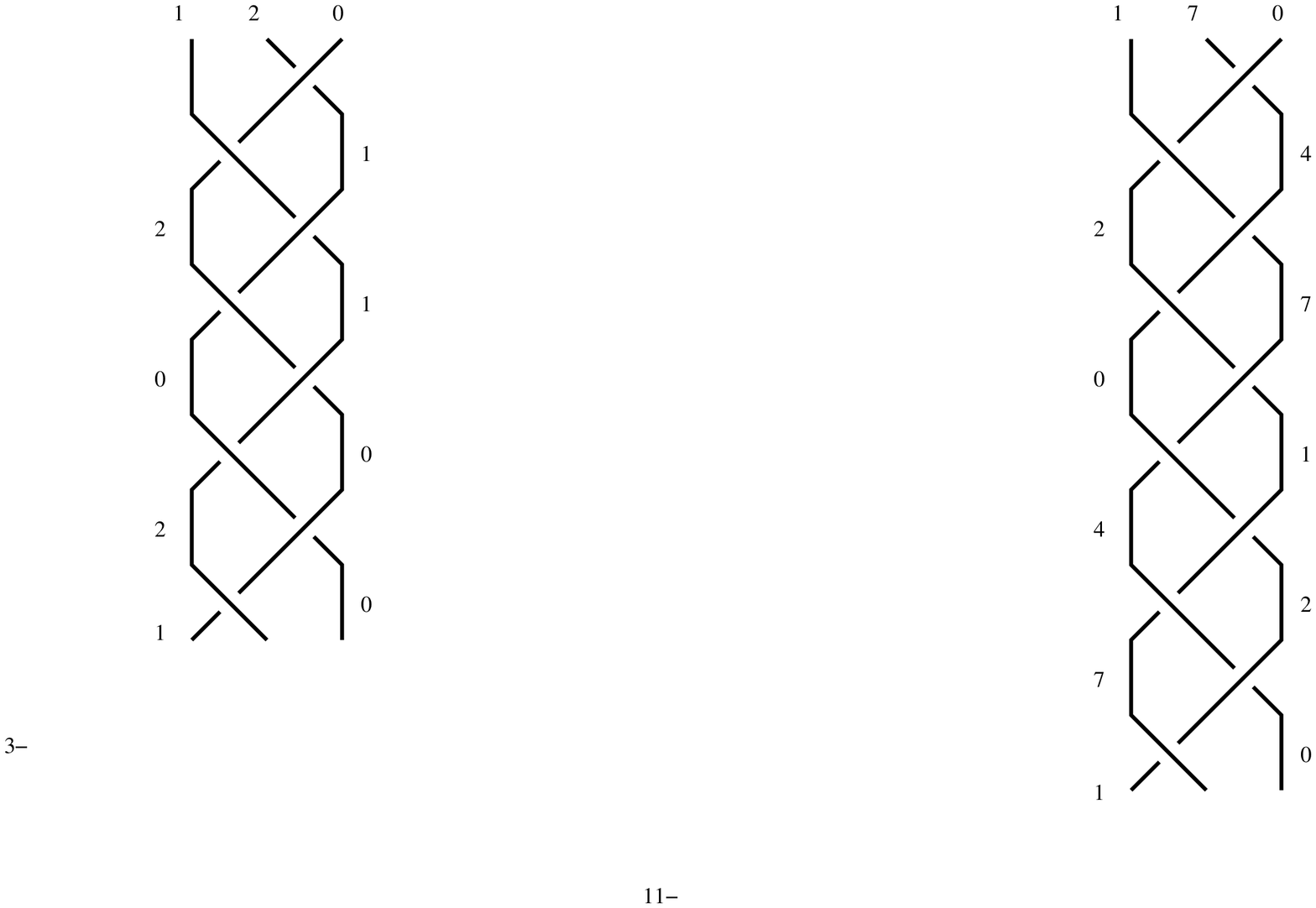}}}
    \caption{Illustrative examples of non-trivial colorings  of $THK(3, n)$'s (cont'd)}\label{fig:thmminexact2}
\end{figure}

\begin{figure}[!ht]
    \psfrag{0}{\huge $0$}
    \psfrag{1}{\huge $1$}
    \psfrag{2}{\huge $2$}
    \psfrag{3}{\huge $3$}
    \psfrag{4}{\huge $4$}
    \psfrag{5}{\huge $5$}
    \psfrag{6}{\huge $6$}
    \psfrag{7}{\huge $7$}
    \psfrag{t}{\large $2$}
    \psfrag{s}{\large $6$}
    \psfrag{7-}{\huge A $7$-coloring of $THK(3, 8)$ with $7$ colors}
    \centerline{\scalebox{.350}{\includegraphics{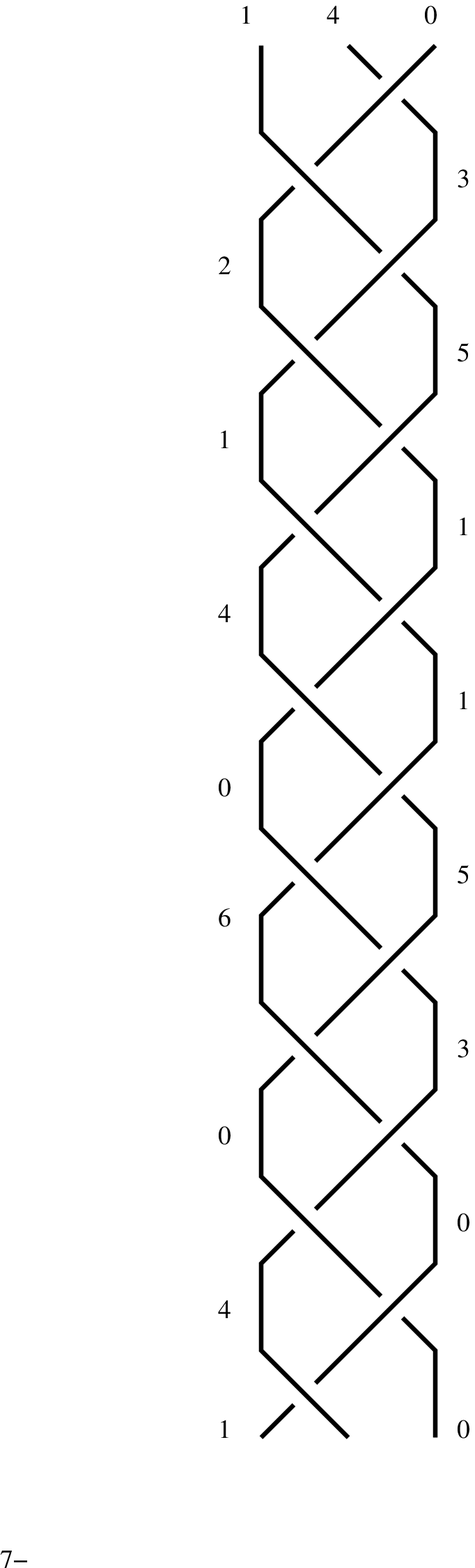}}}
    \caption{Illustrative examples of non-trivial colorings  of $THK(3, n)$'s (end)}\label{fig:thmminexact3}
\end{figure}

After experimenting $r$-colorings of the $THK(3, n)$'s for small values of $r$ and $n$ we came up with the following examples, displayed in Figures \ref{fig:thmminexact1}, \ref{fig:thmminexact2}, and \ref{fig:thmminexact3}. Figures \ref{fig:thmminexact1}, \ref{fig:thmminexact2}, and \ref{fig:thmminexact3} are to be considered upon closure of the braids therein. We do not depict the closure of these braids in order not to over-burden the figures.
Inspection of some cases displayed in Figures \ref{fig:thmminexact1}, \ref{fig:thmminexact2}, and \ref{fig:thmminexact3} shows we have some control on the colors down the right hand-side of the standard diagrams with respect to the left-hand side. This allows us to decrease the number of colors needed to produce a non-trivial coloring.

Consider the $11$-coloring of $THK(3, 5)$ (right-hand side of Figure \ref{fig:thmminexact2}). The sequence of colors down the left-hand side i.e., the $(z\sb{n})$ sequence for $0\leq n\leq 4$, call it $L$, is
\[
L=(1, 2, 0, 4, 7)
\]
The sequence of colors down the middle i.e., the $(y\sb{n})$ sequence for $0\leq n \leq 4$, call it $M$, is
\[
M=(7, 1, 2, 0, 4)
\]
Clearly, $M$ is a circular shift of $L$, since $y\sb{i+1}=x\sb{i}$ for $i=0, 1, 2, 3, 4$, mod $5$. This is a consequence of the arrangement of the arcs in any standard diagram of $THK(3, n)$. Hence, the equality of the sequences $L$ and $M$ modulo circular shift is true for any such diagram.

Let us consider now the sequence of colors down the right-hand side i.e., $(x\sb{n})$ for $0\leq n\leq 4$, call it $R$, of the $11$-coloring of $THK(3, 5)$ (right-hand side of Figure \ref{fig:thmminexact2}):
\[
R = (0, 4, 7, 1, 2)
\]
Then $R$ is a circular shift of $L$ but now this is not a general feature of non-trivial $r$-colorings on standard diagrams of $THK(3,n)$'s (see, for instance, the $5$-coloring of $THK(3, 2)$ on the right-hand side of Figure \ref{fig:thmminexact1}, or the $7$-coloring of $THK(3, 8)$ in Figure \ref{fig:thmminexact3}).

Bearing in mind that if $(r, u\sb{n-1})>1$, there is a non-trivial $r$-coloring of $THK(3, n)$ (Corollary \ref{cor:nontriv}), then, provided $r$ is prime, $\psi (r)$ yields the least number of $\sigma\sb{2}\sigma\sb{1}\sp{-1}$'s we have to juxtapose in order to obtain a non-trivial $r$-coloring for a Turk's head knot on three strands, namely $THK(3, \psi (r))$. For prime $r>5$ and odd $\psi (r)$, we construct a non-trivial $r$-coloring on the standard diagram of $THK(3, \psi (r))$ such that the $R$ sequence is a circular shift of the $L$ sequence (Theorem \ref{thm:minestimates}). There are, on average, two arcs per $\sigma\sb{2}\sigma\sb{1}\sp{-1}$ in a standard diagram of $THK(3, n)$. When $R$ is a circular shift of $L$, we use only one color per $\sigma\sb{2}\sigma\sb{1}\sp{-1}$, on average i.e., we use only $\psi (r)$ colors. This constitutes a reduction in half on the number of colors with respect to the worst case (different arcs, different colors). Moreover, we show that, for prime $r>5$
\[
\psi(r) \mid (r+1) \qquad \text{ or } \qquad \psi(r) \mid (r-1)
\]
so that, if $\psi(r)$ is odd, then $\psi (r)\leq \frac{r+1}{2}$ or $\psi (r)\leq \frac{r-1}{2}$ (Proposition \ref{cor:pmidup}). This guarantees also that we are using roughly half of the colors available, when $r>5$ is prime and $\psi (r)$ is odd.

On the other hand, for prime $r>5$ with $\psi (r)$ even, we show below (Theorem \ref{thm:minestimates}) that the input color $(0, 1, 0)$ induces a coloring whose number of colors is less than $\psi (r)-1$. Here however, it may happen that $\psi (r) = r+1$, so that the estimate equals the total number of colors available.

We develop these ideas below.

\bigbreak

\subsubsection{Further analysis of $u\sb{n}$}\label{subsubsect:unbis}

\noindent

In order to carry out the ideas expressed above, we need a deeper analysis of the sequence $u\sb{n}$. This is the goal of the current Subsection.

\begin{proposition}\label{prop:pmidup}
Let $p\neq 5$ be an odd prime. Then:
\begin{equation*}
\begin{cases}
p\mid u_{p} &\text{if and only if } \quad 5^{\frac{p-1}{2}}=\sb{p} -1\\
p\mid u_{\frac{p-3}{2}} &\text{if and only if } \quad 5^{\frac{p-1}{2}}=\sb{p} 1
\end{cases}
\end{equation*}
\end{proposition}

\begin{proof}
As $u_{n}$ is a sequence taking on integer values we may conclude that the coefficients of $\sqrt{5}$, after applying \textit{Binomial Theorem} to the expression $(\ref{bu})$, will sum zero. Also, one may see that $p\nmid\binom{p+2}{j}$ if and only if $j=0,1,2,p,p+1,p+2$ and $p\nmid\binom{p}{j}$ if and only if $j=0,p$. By \textit{Fermat's Little Theorem} we have $5^{p-1}=\sb{p}1$, which implies $5^{\frac{p-1}{2}}=\sb{p} 1$ or $5^{\frac{p-1}{2}}=\sb{p} -1$. Therefore:

\begin{align*} 0&=\sb{p}u_{p} =\sb{p} \frac{1}{2^{p+2}}\left(2\sum_{i=0}^{\frac{p+1}{2}}\binom{p+2}{2i}5^{\frac{p+1}{2}-i}-
8\sum_{i=0}^{\frac{p-1}{2}}\binom{p}{2i}5^{\frac{p-1}{2}-i}\right) =\sb{p} \\
&=\sb{p} \frac{1}{2^{3}}\left( 2 \times 5^{\frac{p+1}{2}}+2 \times \frac{(p+2)(p+1)}{2}\times 5^{\frac{p-1}{2}}+2(p+2)-8\times 5^{\frac{p-1}{2}}\right)  =\sb{p} \\
&=\sb{p}\frac{1}{2^{3}}\left( 4\times (5^{\frac{p-1}{2}}  +1)\right)
\end{align*}
which is equivalent to$5^{\frac{p-1}{2}}=\sb{p} -1$.

Assume, now, $5^{\frac{p-1}{2}} =\sb{p} 1$. Then by \textit{Euler's Criterion} there exists $\alpha \in \mathbb{Z}$,
such that $\alpha^{2} =\sb{p} 5$. As in $(\ref{bu})$ the coefficients of $\sqrt{5}$ will sum zero, and we obtain:

\begin{align*}
u_{\frac{p-3}{2}}=\sb{p} \alpha^{-1}\left((1+\alpha)^{\frac{p+1}{2}}-4(-1+\alpha)^{\frac{p-3}{2}}-(1-\alpha)^{\frac{p+1}{2}}+4(-1-\alpha)^{\frac{p-3}{2}}\right)
\end{align*}

Multiplying both sides by $2^{\frac{p+1}{2}}\alpha(1+\alpha)^{\frac{p+1}{2}}$, applying \textit{Fermat's Little Theorem} and noting
$p \nmid (\alpha \pm 1)$ for $\alpha^{2} =_{p} 5 \neq 1$:
\begin{align*} 2^{\frac{p+1}{2}}\alpha &(1+\alpha)^{\frac{p+1}{2}} u_{\frac{p-3}{2}}=\sb{p}\\
&=\sb{p}(1+\alpha)^{p+1}-4(1+\alpha)^{2}(-1+\alpha^{2})^{\frac{p-3}{2}}-(1-\alpha^{2})^{\frac{p+1}{2}}+4(-1)^{\frac{p-3}{2}}(1+\alpha)^{p-1}=\sb{p}\\
&=\sb{p} (1+\alpha)^{2}-4(1+\alpha)^{2}(-1+5)^{\frac{p-3}{2}}-(1-5)^{\frac{p+1}{2}}+4(-1)^{\frac{p-3}{2}} =\sb{p} \\
&=\sb{p}  (1+\alpha)^{2}-(1+\alpha)^{2}2^{p-1}-(-1)^{\frac{p+1}{2}}2^{p+1}+4(-1)^{\frac{p-3}{2}} =\sb{p} \\
& =\sb{p}(1+\alpha)^{2}-(1+\alpha)^{2}-4(-1)^{\frac{p-3}{2}}+4(-1)^{\frac{p-3}{2}} =\sb{p} 0
\end{align*}
which yields $u_{\frac{p-3}{2}}=\sb{p}0$, concluding the proof.
\end{proof}

\bigbreak

\begin{corollary}\label{cor:pmidup}
Let $p\neq 5$ be an odd prime. Then:
\begin{equation*}
\begin{cases}
\psi (p)\mid (p+1) &\text{if and only if } \quad 5^{\frac{p-1}{2}}\equiv_{p} -1\\
\psi (p)\mid \big(\frac{p-1}{2}\big) &\text{if and only if } \quad 5^{\frac{p-1}{2}}\equiv_{p} 1
\end{cases}
\end{equation*}
In particular, when $\psi (p)$ is odd, then $\psi (p) \leq (p+1)/2$.
\end{corollary}

\begin{proof} It is a straightforward application of Propositions \ref{proposition:divide} and \ref{prop:pmidup}
\end{proof}

\bigbreak

\begin{proposition}
Let $p$ be a prime greater than $5$.
\begin{enumerate}
\item If $\psi(p)$ is odd then $mincol\sb{p} THK(3, \psi(p)) \leq \psi(p)$
\item If $\psi(p)$ is even then
\begin{enumerate}
\item If $4\mid \psi(p)$, then $mincol\sb{p} THK(3, \psi(p)) \leq \psi(p)-1$
\item If $4\nmid \psi(p)$, then $mincol\sb{p} THK(3, \psi(p)) \leq \psi(p)-5$
\end{enumerate}
\end{enumerate}
\end{proposition}

\begin{proof}
Let us first suppose $\psi(p)$ is odd, and set $a=1$, $c=0$ and $b$ arbitrary, for the moment (Figure \ref{fig:slqbis}).
\begin{figure}[!ht]
    \psfrag{a}{\huge $1$}
    \psfrag{b}{\huge $b$}
    \psfrag{c}{\huge $0$}
    \psfrag{x}{\huge $x_1 = 2$}
    \psfrag{y}{\huge $y_1=1$}
    \psfrag{z}{\huge $z_1 = -b$}
    \centerline{\scalebox{.50}{\includegraphics{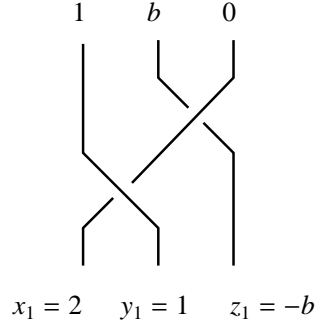}}}
    \caption{Propagation of colors down $\sigma\sb{2}\sigma\sb{1}\sp{-1}$ for a particular choice of $a, b, c$}\label{fig:slqbis}
\end{figure}
We will show that $b$ can be chosen so that the $R$ sequence i.e., the sequence $(z\sb{n})|\sb{0\leq n\leq \psi (p)}$, is a circular shift of the $L$ sequence i.e., the sequence $(x\sb{n})|\sb{0\leq n\leq \psi (p)}$, yielding the required result. As we saw above, in the beginning of this Subsection, this implies that the number of colors is bounded above by $\psi (p)$.

By Proposition \ref{proposition:const}:
\begin{equation}\label{eqn : eqn}
x_{k}=3x_{k-1}-x_{k-2}-1+b \qquad \Leftrightarrow \qquad (x_{k}-1+b)=3(x_{k-1}-1+b)-(x_{k-2}-1+b)
\end{equation}
and similarly for $y_{k},z_{k}$. In particular, the above relations tell us that any term of $x\sb{k}$ (respect., $z\sb{k}$) is obtained from the preceding two terms of the sequence. We then focus on two consecutive terms of the $x$-sequence, $(x\sb{k}, x\sb{k+1})$ and we will later equate a certain pair of them to $(z\sb{0}, z\sb{1})=(0, -b)$ in order to obtain the equality between the $R$ and $L$ sequences, modulo circular shift.

The expressions in \ref{eqn : eqn} lead to:
\begin{align*}
\left[
\begin{matrix}
x_{k+1}\\
x_{k}
\end{matrix}
\right]&= \left[
\begin{matrix}
3 &-1\\
1 &0
\end{matrix}
\right]
\left[
\begin{matrix}
x_{k}-1+b\\
x_{k-1}-1+b
\end{matrix}
\right]+\left[
\begin{matrix}
1-b\\
1-b
\end{matrix}
\right]= \left[
\begin{matrix}
3 &-1\\
1 &0
\end{matrix}
\right]
\left( \left[
\begin{matrix}
3 &-1\\
1 &0
\end{matrix}
\right]
\left[
\begin{matrix}
x_{k-1}-1+b\\
x_{k-2}-1+b
\end{matrix}
\right]\right) +\left[
\begin{matrix}
1-b\\
1-b
\end{matrix}
\right]\\
&= \dots =\left[
\begin{matrix}
3 &-1\\
1 &0
\end{matrix}
\right]^{k}
\left[
\begin{matrix}
x_{1}-1+b\\
x_{0}-1+b
\end{matrix}
\right]+\left[
\begin{matrix}
1-b\\
1-b
\end{matrix}
\right]=
\left[
\begin{matrix}
3 &-1\\
1 &0
\end{matrix}
\right]^{k}
\left[
\begin{matrix}
1+b\\
b
\end{matrix}
\right]+\left[
\begin{matrix}
1-b\\
1-b
\end{matrix}
\right]
\end{align*}

By induction one can prove that:

\begin{equation*}
\left[
\begin{matrix}
3 &-1\\
1 &0
\end{matrix}
\right]^{k}=
\left[
\begin{matrix}
u_{2k+1} &-u_{2k-1}\\
u_{2k-1} &-u_{2k-3}
\end{matrix}
\right]
\end{equation*}

and taking determinants on both sides of the preceding equation we obtain:

\begin{equation*}
1 = -u_{2k+1}u_{2k-3}+u_{2k-1}^{2}
\end{equation*}

In order for $L$ to be a circular shift of $R$ we set $x_{k+1}=z_{1}=-b,x_{k}=z_{0}=0$. Then:
\begin{align*}
\left[
\begin{matrix}
-b\\
0
\end{matrix}
\right]=
\left[
\begin{matrix}
3 &-1\\
1 &0
\end{matrix}
\right]^{k}
\left[
\begin{matrix}
1+b\\
b
\end{matrix}
\right]+\left[
\begin{matrix}
1-b\\
1-b
\end{matrix}
\right]
\end{align*}
which is equivalent to:
\begin{equation}\label{eq:b}
\left[
\begin{matrix}
u_{2k+1}+1 &-u_{2k-1}-1\\
u_{2k-1}+1 &-u_{2k-3}-2
\end{matrix}
\right]
\left[
\begin{matrix}
1+b\\
b
\end{matrix}
\right]\equiv_{p}
\left[
\begin{matrix}
0\\
0
\end{matrix}
\right]
\end{equation}

The determinant of the coefficient matrix in (\ref{eq:b}) is given by:
\begin{multline*}
 -u_{2k+1}u_{2k-3}+u_{2k-1}^{2}-2u_{2k+1}-u_{2k-3}+2u_{2k-1}-1=\\
=1+(-u_{2k-3}+3u_{2k-1}-u_{2k-1})-u_{2k+1}-u_{2k-1}-1=\\
=-(u_{2k+1}+u_{2k-1})=-u_{2k}
\end{multline*}
We then set $k=\frac{\psi(p)-1}{2}$ so that $2k+1=\psi (p)$. Then $u\sb{2k}=\sb{p}0$, by definition of $\psi $, and the kernel associated to the system \ref{eq:b} is non-trivial. If the non-null vectors of this kernel had equal coordinates then in particular,
\begin{equation*}
\left[
\begin{matrix}
0\\
0
\end{matrix}
\right]=
\left[
\begin{matrix}
u_{\psi(p)}+1 &-u_{\psi(p)-2}-1\\
* &*
\end{matrix}
\right]
\left[
\begin{matrix}
1\\
1
\end{matrix}
\right]=\left[
\begin{matrix}
u_{\psi(p)}-u_{\psi(p)-2}\\
*
\end{matrix}
\right]= \dots
\end{equation*}
and, by Proposition \ref{proposition:sominhas}, as $\psi(p)$ is odd:
\begin{equation*}
\dots =\left[
\begin{matrix}
u_{\psi(p)}-u_{\psi(p)-2}-(u_{\psi(p)}-u_{\psi(p)-1}+u_{\psi(p)-2})\\
*
\end{matrix}
\right]=_{p}
\left[
\begin{matrix}
-2u_{\psi(p)-2}\\
*
\end{matrix}
\right]
\end{equation*}
Now, $2u_{\psi(p)-2}=_{p} 0$, only if $p=2$ (in which case the proposition is easily verified, as $\psi(2)=3\geq 2$), or $p\mid u_{\psi(p)-2}$, following from Proposition \ref{proposition:divide}  that $\psi(p) \mid \psi(p)-1$, which implies $\psi(p)=1$ and $p\mid u_{0}=1$, contradicting the fact that $p$ is prime. We can thus assume that non-trivial elements in the indicated kernel have distinct coordinates, say $(r, s)$ with $r\neq s$ mod $p$.

The vector $((r-s)^{-1}r,(r-s)^{-1}s)=_{p}(r',s')$ also belongs to this kernel and satisfies $r'-s'=1$. Therefore
$b=s'$ is a solution to equation (\ref{eq:b}). We conclude that, for the top colors $(1,s',0)$, the rightmost and leftmost arcs have equal pairs of consecutive colors. Therefore, given the recurrence relation satisfied by the colors on each side, the colors used in the leftmost, middle, and rightmost arcs are the same and we use, at most, $\psi(p)$ colors.

\bigbreak

Let us now suppose that $\psi(p)$ is even. By Proposition \ref{proposition:sominhas} we have:

\begin{equation*}
\frac{u_{\psi(p)}+u_{\psi(p)-2}}{5}=u_{\psi(p)-1} =_{p}0  \quad \Rightarrow  \quad u_{\psi(p)}=_{p} -u_{\psi(p)-2}
\end{equation*}
as $p\neq 5$. Remember that $u_{0}=1=-u_{-2}$. So, for some integer $\alpha$ we have $u_{\psi(p)-2}=_{p} \alpha u_{0}$
and $u_{\psi(p)}=_{p} \alpha u_{-2}$.
Furthermore, and as we have:
\begin{equation*}
u_{k}=3u_{k-2}-u_{k-4} \quad  \Leftrightarrow \quad u_{k-4}=3u_{k-2}-u_{k}
\end{equation*}
one can easily see that,
\begin{align*}
\left[
\begin{matrix}
u_{-2}\\
u_{0}
\end{matrix}
\right]&=
\left[
\begin{matrix}
3 &-1\\
1 &0
\end{matrix}
\right]
\left[
\begin{matrix}
u_{0}\\
u_{2}
\end{matrix}
\right]
=
\left[
\begin{matrix}
3 &-1\\
1 &0
\end{matrix}
\right]
\left[
\begin{matrix}
3 &-1\\
1 &0
\end{matrix}
\right]
\left[
\begin{matrix}
u_{0}\\
u_{2}
\end{matrix}
\right]
= \cdots =
\left[
\begin{matrix}
3 &-1\\
1 &0
\end{matrix}
\right]^{\frac{\psi(p)}{2}}
\left[
\begin{matrix}
u_{\psi(p)-2}\\
u_{\psi(p)}
\end{matrix}
\right] =\\
&=_{p}
\alpha\left[
\begin{matrix}
3 &-1\\
1 &0
\end{matrix}
\right]^{\frac{\psi(p)}{2}}
\left[
\begin{matrix}
u_{0}\\
u_{-2}
\end{matrix}
\right]=
\alpha\left[
\begin{matrix}
3 &-1\\
1 &0
\end{matrix}
\right]^{\frac{\psi(p)}{2}-1}
\left[
\begin{matrix}
u_{2}\\
u_{0}
\end{matrix}
\right]= \cdots =
\alpha
\left[
\begin{matrix}
u_{\psi(p)}\\
u_{\psi(p)-2}
\end{matrix}
\right]
=_{p}
\alpha^{2}
\left[
\begin{matrix}
u_{-2}\\
u_{0}
\end{matrix}
\right]
\end{align*}

Then $\alpha^{2}=_{p} 1\Leftrightarrow \alpha=_{p} \pm 1$. We will now show that $\alpha =\sb{p}1$.
\bigbreak

Assume to the contrary and suppose $\alpha =\sb{p}-1$. Then
$u_{0}=_{p} -u_{\psi(p)-2}$, $u_{2}=_{p} 3u_0-u\sb{-2}=\sb{p}-3u\sb{\psi (p)-2}+u\sb{\psi (p)}=\sb{p}-u_{\psi(p)-4}$ and, by applying Proposition \ref{proposition:sominhas} twice,
we get $u_{1}=1=_{p} -u_{\psi(p)-3}$. Using the recurrence satisfied by $u_{n}$, more generally, we have,
\begin{equation*}
u_{-2+i}=_{p} -u_{\psi(p)-i}\text{ , }i=0,1,...,\psi(p)+2
\end{equation*}
In particular, $\alpha =_{p}-1$ implies $u_{-2+\frac{\psi(p)}{2}+1}=_{p}-u_{\psi(p)-\frac{\psi(p)}{2}-1}
 \Leftrightarrow u_{\frac{\psi(p)}{2}-1}=_{p} 0$, which contradicts  the definition of $\psi(p)$.

\bigbreak

Thus $\alpha = 1$ and $u_{\psi(p)-2}=_{p} 1 =_{p} -u_{\psi(p)}$.

\bigbreak


We note that $(x\sb{0}, y\sb{0}, z\sb{0})=(0, 1, 0)$ induces a non-trivial $p$-coloring on $THK(3, \psi (p))$, since, by definition, $p\mid \psi u\sb{(p)-1}$. Then, $(x\sb{\psi (p)}, y\sb{\psi (p)}, z\sb{\psi (p)})=(0, 1, 0)$.

Arguing as above in the odd $\psi (p)$ case, we obtain:
\begin{align*}
\left[
\begin{matrix}
x_{k+1}\\
x_{k}
\end{matrix}
\right]&=
\left[
\begin{matrix}
3 &-1\\
1 &0
\end{matrix}
\right]^{k}
\left[
\begin{matrix}
x_{1}+1\\
x_{0}+1
\end{matrix}
\right]+\left[
\begin{matrix}
-1\\
-1
\end{matrix}
\right]=
\left[
\begin{matrix}
u_{2k+1} &-u_{2k-1}\\
u_{2k-1} &-u_{2k-3}
\end{matrix}
\right]
\left[
\begin{matrix}
1\\
1
\end{matrix}
\right]+\left[
\begin{matrix}
-1\\
-1
\end{matrix}
\right]\\
\left[
\begin{matrix}
z_{k+1}\\
z_{k}
\end{matrix}
\right]&=
\left[
\begin{matrix}
3 &-1\\
1 &0
\end{matrix}
\right]^{k}
\left[
\begin{matrix}
z_{1}+1\\
z_{0}+1
\end{matrix}
\right]+\left[
\begin{matrix}
-1\\
-1
\end{matrix}
\right]=
\left[
\begin{matrix}
u_{2k+1} &-u_{2k-1}\\
u_{2k-1} &-u_{2k-3}
\end{matrix}
\right]
\left[
\begin{matrix}
0\\
1
\end{matrix}
\right]+\left[
\begin{matrix}
-1\\
-1
\end{matrix}
\right]
\end{align*}

Using these expressions along with $u\sb{\psi (p)-2}=\sb{p}1=\sb{p}-u\sb{\psi (p)}$ and Proposition \ref{proposition:sominhas}, we obtain the colors displayed in Figure \ref{fig:psipeven}.
\begin{figure}[!ht]
    \psfrag{y0}{\huge $1$}
    \psfrag{x0}{\huge $x\sb{0}=0$}
    \psfrag{x1}{\huge $x\sb{1}=0$}
    \psfrag{x2}{\huge $x\sb{2}=1$}
    \psfrag{z0}{\huge $z\sb{0}=0$}
    \psfrag{z1}{\huge $z\sb{1}=-1$}
    \psfrag{z2}{\huge $z\sb{2}=-2$}
    \psfrag{xpsi2}{\huge $x\sb{\frac{\psi (p)}{2}}=-2$}
    \psfrag{xpsi21}{\huge $x\sb{\frac{\psi (p)}{2}+1}=-2$}
    \psfrag{zpsi2}{\huge $z\sb{\frac{\psi (p)}{2}}=-2$}
    \psfrag{zpsi21}{\huge $z\sb{\frac{\psi (p)}{2}+1}=-1$}
    \psfrag{zpsi22}{\huge $z\sb{\frac{\psi (p)}{2}+2}=0$}
    \psfrag{zpsi23}{\huge $z\sb{\frac{\psi (p)}{2}+3}=2$}
    \psfrag{xpsi11}{\huge $x\sb{\psi (p)-1}=1$}
    \psfrag{xpsi1}{\huge $x\sb{\psi (p)}=0$}
    \psfrag{zpsi11}{\huge $z\sb{\psi (p)-1}=2$}
    \psfrag{zpsi1}{\huge $z\sb{\psi (p)}=0$}
    \psfrag{-3}{\huge $-3$}
    \psfrag{-6}{\huge $-6$}
    \psfrag{...}{\huge $\vdots$}
    \centerline{\scalebox{.50}{\includegraphics{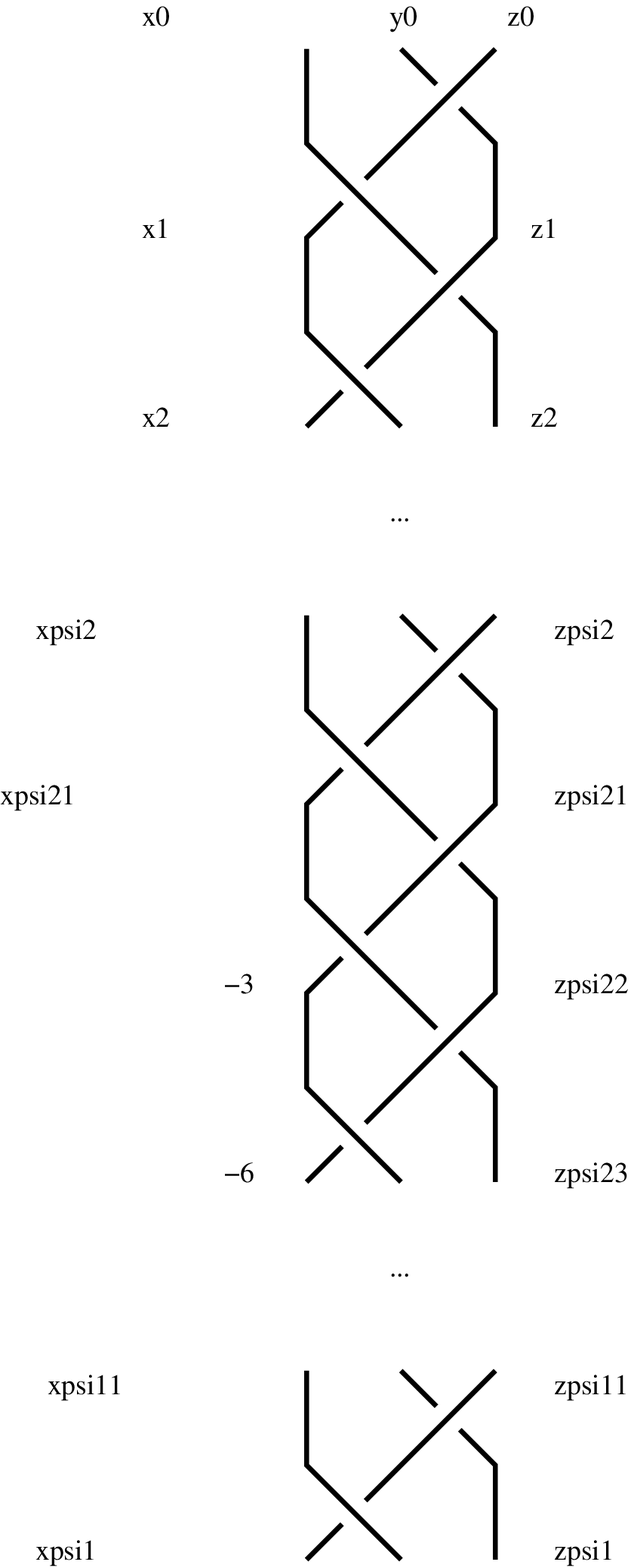}}}
    \caption{Minimizing the number of colors for prime $p$ with $\psi (p)$ even.}\label{fig:psipeven}
\end{figure}

By Proposition \ref{proposition:const}, the following hold:
\begin{align*}
x_{k}&=3x_{k-1}-x_{k-2}+1 \quad \Leftrightarrow  \quad x_{k-2}=3x_{k-1}-x_{k}+1\\
z_{k}&=3z_{k-1}-z_{k-2}+1  \quad \Leftrightarrow  \quad z_{k-2}=3z_{k-1}-z_{k}+1
\end{align*}

Since $x_{1} = 0 = x_{\psi(p)}$ and $x\sb{2}=1=x_{\psi(p)-1}$, then, the sequence $(x_{i})\mid\sb{1\leq i\leq \psi(p)}$ contains at most $\frac{\psi(p)}{2}$ distinct terms.

\bigbreak

Assume $4\mid \psi (p)$ and consider the sequence $(z_{i})\mid\sb{1\leq i\leq \frac{\psi(p)}{2}+1}$. We have, $z_{1} = -1 = z_{\frac{\psi(p)}{2}+1}$ and $z\sb{2}=-2=z_{\frac{\psi(p)}{2}}$, then there are equal terms in two's. Since $\frac{\psi(p)}{2}+1$ is odd, this means there is a central term, $z\sb{\frac{\psi(p)}{4}+1}$, and the remaining terms are equal in two's. Thus, this sequence contributes with $\frac{\psi(p)}{4}+1$ distinct terms.

Now for the sequence $(z_{i})\mid\sb{\frac{\psi(p)}{2}+2\leq i\leq \psi(p)}$. Reasoning as in the preceding case, this sequence has at most, $\frac{\psi(p)-(\psi(p)/2+2)}{2}+1=\frac{\psi (p)}{4}$.

So, in this $4\mid \psi (p)$ instance, the number of distinct colors is at most
\[
\frac{\psi(p)}{2} + \left( \frac{\psi(p)}{4}+1\right) + \frac{\psi (p)}{4} - 2 = \psi (p) -1
\]
where the $-2$ stems from the fact that $0$ and $-2$ appear on both the $x$ and the $z$ sequences (see Figure \ref{fig:psipeven}).







Assume now $4\nmid \psi (p)$ and consider the sequence $(z_{i})\mid\sb{1\leq i\leq \frac{\psi (p)}{2}+1}$.  Here $\frac{\psi (p)}{2}+1$ is even, so this sequence is made up of $\frac{\psi(p)/2+1}{2}$ pairs of equal terms. Moreover, there are two equal consecutive terms,  $z\sb{\frac{\psi(p)/2+1}{2}}=z\sb{\frac{\psi(p)/2+3}{2}}$ (see Figure \ref{fig:psipeveneven}, with $c$ the common value). Given the arrangement of the arcs in these standard diagrams, this implies that the $z\sb{\frac{\psi(p)/2+1}{2}}=\sb{p}c$ (see Figure \ref{fig:psipeveneven}). Using the coloring condition at the crossings and the recurrence relation $z\sb{k-2}=3z\sb{k-1}-z\sb{k}+1$, we obtain the other colors in Figure \ref{fig:psipeveneven}, which tell us that there are two colors in $(z_{i})\mid\sb{1\leq i\leq \frac{\psi (p)}{2}+1}$ that already showed up in the $x$ sequence: $c$ and $-1$. Thus, the net contribution of $(z_{i})\mid\sb{1\leq i\leq \frac{\psi (p)}{2}+1}$ is at most $\frac{\psi(p)/2+1}{2}-2$.

Now, for the sequence $(z_{i})\mid\sb{\frac{\psi (p)}{2}+2 \leq i\leq \psi (p)}$. Again, we have an even number of terms, $\frac{\psi (p)}{2}-1$, equal in two's. Since $z\sb{\frac{\psi (p)}{2}+2}=\sb{p}0=\sb{p}z\sb{\psi (p)}$, which has already been accounted for in the $x$ sequence, then we will consider only $\frac{\frac{\psi (p)}{2}-3}{2}$ pairs of terms. Moreover, there is again a phenomenon analogous to the one depicted in Figure \ref{fig:psipeveneven}, induced by the fact that there are two equal consecutive terms in the $(z_{i})\mid\sb{\frac{\psi (p)}{2}+2 \leq i\leq \psi (p)}$ sequence. This time is just the $c$ that has to be considered. Moreover, we will now prove that this $c$ is distinct from the $c$ occurring in connection with the sequence $(z_{i})\mid\sb{1\leq i\leq \frac{\psi (p)}{2}+1}$. If they were equal, since the next term in both subsequences, $(x_{i})\mid\sb{1\leq i\leq \frac{\psi(p)}{2}}$ and $(x_{i})\mid\sb{\frac{\psi(p)}{2}+1\leq i\leq \psi(p)}$, is $-1$ (mod $p$), then, given the recurrence relation $z\sb{k}=3z\sb{k-1}-z\sb{k-2}+1$, then the subsequences would be equal. But they are not equal. Hence the net contribution of $(z_{i})\mid\sb{\frac{\psi (p)}{2}+2 \leq i\leq \psi (p)}$ to the total number of colors is at most $\frac{\frac{\psi (p)}{2}-3}{2}-1$.

Finally, the total number of colors in this $4\nmid \psi (p)$ instance is at most (recall that $-2$ shows up on both sides):
\[
\frac{\psi(p)}{2} + \left( \frac{\psi(p)/2+1}{2}-2\right) + \left( \frac{\frac{\psi (p)}{2}-3}{2}-1\right) -1 = \psi (p) -5
\]

\begin{figure}[!ht]
    \psfrag{xpsi21}{\huge $c=x\sb{\frac{\psi (p)-1}{4}-1}$}
    \psfrag{zpsi21}{\huge $z\sb{\frac{\psi (p)-1}{4}-1}=3c-c+1=2c+1$}
    \psfrag{zpsi22}{\huge $z\sb{\frac{\psi (p)-1}{4}}=c$}
    \psfrag{zpsi23}{\huge $z\sb{\frac{\psi (p)-1}{4}+1}=c$}
    \psfrag{-3}{\huge $-1=2c-(2c+1)=x\sb{\frac{\psi (p)-1}{4}-1}$}
    \psfrag{...}{\huge $\vdots$}
    \centerline{\scalebox{.50}{\includegraphics{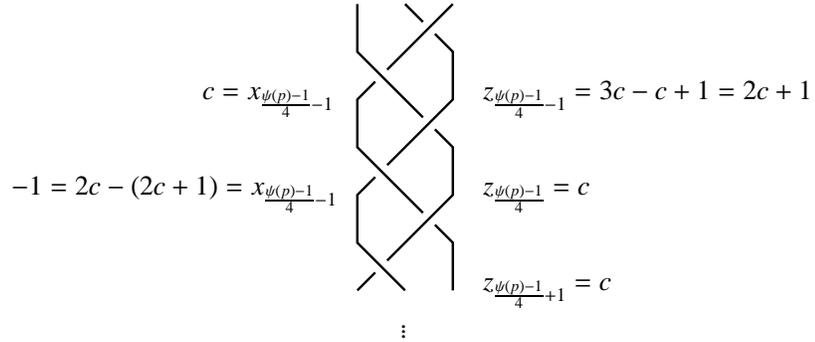}}}
    \caption{The colors around $z\sb{\frac{\psi(p)}{2}+1}, z\sb{\frac{\psi(p)}{2}+1+1}$ for prime $p$ with $\psi (p)$ even and $4\nmid \psi (p)$.}\label{fig:psipeveneven}
\end{figure}

This concludes the proof for the $\psi(p)$ even case, concluding the proof of Proposition \ref{prop:pmidup}.

\end{proof}

\bigbreak

\subsubsection{Proof of Theorem \ref{thm:minestimates} and Corollary \ref{cor:minestimatesgen}}\label{subsubsect:thmminestimates}

\noindent

The proof of Theorem \ref{thm:minestimates} is a straightforward application of Proposition \ref{prop:pmidup}.

As for proof of Corollary \ref{cor:minestimatesgen}, given positive integers $n$ and $r>1$, we choose a common prime factor of $r$ and $u\sb{n-1}$ (in order to ensure there are non-trivial colorings), which minimizes $\psi$ (in order to minimize the number of colors involved. Applying Propositions \ref{prop:stackcol}  and \ref{prop:pcolrcol}, we conclude the proof of Corollary \ref{cor:minestimatesgen}.

\bigbreak

\clearpage

\appendix

\section{Table of $\psi$}\label{app:tablepsi}
\noindent

\begin{table}[h!!]
\begin{center}
\begin{tabular}{| l | r |||||||| l | r |||||||| l | r |||||||| l | r |||||||| l | r |}
\hline
$\mathbf{r}$ & $\mathbf{\psi (r)}$ &   $\mathbf{r}$ & $\mathbf{\psi (r)}$ &   $\mathbf{r}$ & $\mathbf{\psi (r)}$  &   $\mathbf{r}$ & $\mathbf{\psi (r)}$ &   $\mathbf{r}$ & $\mathbf{\psi (r)}$ \\ \hline\hline
$1$  & $2$ &  $38$  & $9$  &  $75$ &    $100$  &  $112$ &    $24$ &  $\mathbf{149}$ &    $\mathbf{74}$ \\ \hline
$\mathbf{2}$  & $\mathbf{3}$  &  $39$  & $28$ &  $76$  & $9$  &  $\mathbf{113}$ &    $\mathbf{38}$ &  $150$ &    $300$ \\ \hline
$\mathbf{3}$  & $\mathbf{4}$  &  $40$  & $30$ &  $77$  & $40$  &  $114$ &    $36$ &  $\mathbf{151}$ &    $\mathbf{25}$ \\ \hline
$4$  & $3$  &  $\mathbf{41}$  & $\mathbf{20}$  &  $78$  & $84$  &  $115$ &    $120$ &  $152$ &    $18$ \\ \hline
$\mathbf{5}$  & $\mathbf{10}$  &  $42$ & $24$  &  $\mathbf{79}$  & $\mathbf{39}$  &  $116$ &    $21$ &  $153$ &    $36$ \\ \hline
$6$  & $12$  &  $\mathbf{43}$  & $\mathbf{44}$  &  $80$  & $60$  &  $117$ &    $84$ &  $154$ &    $120$ \\ \hline
$\mathbf{7}$  & $\mathbf{8}$  &  $44$  & $15$  &  $81$  & $108$  &  $118$ &    $87$&  $155$ &    $30$ \\ \hline
$8$  & $6$  &  $45$  & $60$  &  $82$  & $60$  &  $119$ &    $72$ &  $156$ &    $84$ \\ \hline
$9$  & $12$  &  $46$  & $24$  &  $\mathbf{83}$  & $\mathbf{84}$  &  $120$ &    $60$ &  $\mathbf{157}$ &    $\mathbf{158}$ \\ \hline
$10$  & $30$  &  $\mathbf{47}$  & $\mathbf{16}$  &  $84$  & $24$  &  $121$ &    $55$ &  $158$ &    $39$ \\ \hline
$\mathbf{11}$  & $\mathbf{5}$  &  $48$  & $12$  &  $85$  & $90$  &  $122$ &    $30$ &  $159$ &    $108$ \\ \hline
$12$  & $12$  &  $49$  & $56$  &  $86$  & $132$  &  $123$ &    $20$ &  $160$ &    $120$ \\ \hline
$\mathbf{13}$  & $\mathbf{14}$  &  $50$  & $150$  &  $87$  & $28$  &  $124$ &    $15$ &  $161$ &    $24$ \\ \hline
$14$  & $24$  &  $51$  & $36$  &  $88$  & $30$  &  $125$ &    $250$ &  $162$ &    $28$ \\ \hline
$15$  & $20$  &  $52$  & $42$  &  $\mathbf{89}$  & $\mathbf{22}$  &  $126$ &    $24$ &  $\mathbf{163}$ &    $\mathbf{164}$ \\ \hline
$16$  & $12$  &  $\mathbf{53}$  & $\mathbf{54}$  &  $90$  & $60$  &  $\mathbf{127}$ &    $\mathbf{128}$ &  $164$ &    $60$ \\ \hline
$\mathbf{17}$  & $\mathbf{18}$  &  $54$  & $36$  &  $91$  & $56$  &  $128$ &    $96$ &  $165$ &    $20$ \\ \hline
$18$  & $12$  &  $55$  & $10$  &  $92$  & $24$  &  $129$ &    $44$ &  $166$ &    $84$ \\ \hline
$\mathbf{19}$  & $\mathbf{9}$  &  $56$  & $24$  &  $93$  & $60$  &  $130$ &    $210$ &  $\mathbf{167}$ &    $\mathbf{168}$ \\ \hline
$20$  & $30$  &  $57$  & $36$  &  $94$  & $48$  &  $\mathbf{131}$ &    $\mathbf{65}$ &  $168$ &    $24$ \\ \hline
$21$  & $8$  & $58$  & $21$ &  $95$  & $90$   &  $132$ &    $60$ &  $169$ &    $182$ \\ \hline
$22$  & $15$ & $\mathbf{59}$ & $\mathbf{29}$  &  $96$  & $24$  &  $133$ &    $72$ &  $170$ &    $90$ \\ \hline
$\mathbf{23}$  & $\mathbf{24}$ & $60$ & $60$  &  $\mathbf{97}$  & $\mathbf{98}$  &  $134$ &    $204$ &  $171$ &    $36$ \\ \hline
$24$ & $12$  & $\mathbf{61}$ & $\mathbf{30}$  &  $98$  & $168$  &  $135$ &    $180$ &  $172$ &    $132$ \\ \hline
$25$  & $50$  & $62$ & $15$ &  $99$  & $60$   &  $136$ &    $18$ &  $\mathbf{173}$ &    $\mathbf{174}$ \\ \hline
$26$  & $42$  & $63$ & $24$ &  $100$  & $150$   &  $\mathbf{137}$ &    $\mathbf{138}$ &  $174$ &    $84$ \\ \hline
$27$  & $36$  & $64$ & $48$ &  $\mathbf{101}$  & $\mathbf{25}$   &  $138$ &    $24$ &  $175$ &    $200$ \\ \hline
$28$  & $24$  & $65$ & $70$  &  $102$  & $36$  &  $\mathbf{139}$ &    $\mathbf{23}$ &  $176$ &    $60$ \\ \hline
$\mathbf{29}$  & $\mathbf{7}$  & $66$ & $60$  &  $\mathbf{103}$  & $\mathbf{104}$  &  $140$ &    $120$ &  $177$ &    $116$ \\ \hline
$30$  & $60$  & $\mathbf{67}$ & $\mathbf{68}$  &  $104$  & $42$  &  $141$ &    $16$ &  $178$ &    $66$ \\ \hline
$\mathbf{31}$  & $\mathbf{15}$  & $68$ & $18$  &  $105$  & $40$  &  $142$ &    $105$ &  $\mathbf{179}$ &    $\mathbf{89}$ \\ \hline
$32$  & $24$  & $69$ & $24$ &  $106$  & $54$   &  $143$ &    $70$ &  $180$ &    $60$ \\ \hline
$33$  & $20$  & $70$ & $120$ &  $\mathbf{107}$  & $\mathbf{36}$   &  $144$ &    $12$ &  $\mathbf{181}$ &    $\mathbf{45}$ \\ \hline
$34$  & $18$  & $\mathbf{71}$ & $\mathbf{35}$  & $108$  &  $36$   &  $145$ &    $70$ &  $182$ &    $168$ \\ \hline
$35$  & $40$  & $72$ &    $12$  & $\mathbf{109}$  &  $\mathbf{54}$   &  $146$ &    $222$ &  $183$ &    $60$ \\ \hline
$36$  & $12$ & $\mathbf{73}$ &    $\mathbf{74}$   &  $110$  & $30$   &  $147$ &    $56$ &  $184$ &    $24$ \\ \hline
$\mathbf{37}$  & $\mathbf{38}$ & $74$ &   $114$  & $111$  &  $76$   &  $148$ &    $114$ &  $185$ &    $190$ \\ \hline
\end{tabular}
\caption{Table of $\psi$. In bold: prime $r$'s and their $\psi$'s.}\label{Ta:psi}
\end{center}
\end{table}
\end{document}